\documentclass[11pt]{amsart}

\usepackage{amsmath}
\usepackage{amsfonts}
\usepackage{amsthm}

\usepackage{mathtools}
\mathtoolsset{showonlyrefs}

\usepackage{graphicx} % Required for inserting images
\usepackage{float}
\usepackage{svg}

\usepackage[left=2.65cm,right=2.65cm,top=2.7cm,bottom=2.7cm]{geometry}

\newcommand{\cD}{{\mathfrak{D}}}

\usepackage{xcolor}

\usepackage{hyperref}
\hypersetup{
	colorlinks=true,
	linkcolor=blue,
	citecolor=red,
	filecolor=black,
	urlcolor=black,
}

\usepackage{csquotes}

\newtheorem{theorem}{Theorem}[section]

\newtheorem{proposition}[theorem]{Proposition}
\newtheorem{lemma}[theorem]{Lemma}

\theoremstyle{remark}
\newtheorem{remark}[theorem]{Remark}
\newtheorem{remarks}[theorem]{Remarks}

\theoremstyle{definition}
\newtheorem*{definition*}{Definition}

\numberwithin{equation}{section}

\title[Splitting and Merging of Stagnation Points]{Splitting and Merging of Stagnation Points\\of Solutions to the 2D Navier-Stokes Equations}

\author{Isidro Benaroya} \address{Instituto de Ciencias Matem\'aticas, Consejo Superior de Investigaciones Cient\'\i ficas, 28049 Madrid, Spain}
\email{isidro.benaroya@icmat.es}

\author{Alberto Enciso} \address{Instituto de Ciencias Matem\'aticas, Consejo Superior de Investigaciones Cient\'\i ficas, 28049 Madrid, Spain}
\email{aenciso@icmat.es}

\author{Daniel Peralta-Salas} \address{Instituto de Ciencias Matem\'aticas, Consejo Superior de Investigaciones Cient\'\i ficas, 28049 Madrid, Spain}
\email{dperalta@icmat.es}

\newcommand\bu{\mathbf{u}}
\newcommand\bx{\mathbf{x}}
\DeclareMathOperator{\Div}{div}
\DeclareMathOperator{\rank}{rank}

\begin{document}

\begin{abstract}
We construct solutions to the Navier--Stokes equations on $\mathbb{R}^2$ and $\mathbb{T}^2$ that exhibit an arbitrary number of stagnation points which merge and split along trajectories that can be prescribed freely, up to a small deformation.
\end{abstract}

\maketitle

\setcounter{tocdepth}{1}
\tableofcontents

\section{Introduction}\label{section:intro}

Let us consider the 2D Navier--Stokes equations,
\[
\frac{\partial \bu}{\partial t} + (\bu\cdot\nabla)\bu = -\nabla p +  \Delta \bu, \qquad
\Div  \bu = 0,
\]
which model a viscous incompressible fluid on the plane. Note that the kinematic viscosity has been normalized to~1 and we do not consider any external forces. We can assume that the spatial variable $\bx=(x,y)$ takes values in~$\mathbb R^2$, although later on we will also consider the torus $\mathbb{T}^2:=(\mathbb{R}/2\pi\mathbb{Z})^2$.

In this paper we are interested in the {\em stagnation points} of the fluid, that is, the points where the velocity field $\bu(\mathbf x,t)$ vanishes at time~$t$. Specifically, we aim to study how stagnation points can split and merge as the fluid evolves. 
This topic has been previously studied by various authors, notably Li and Sinai~\cite{li_bifurcations_2012,li_nonsymmetric_2012} and Moffatt~\cite{Mo,moffatt_topological_2021}. In short, they construct solutions to the Navier--Stokes equations on $\mathbb T^2$ whose stream functions exhibit some critical points that merge or split in a certain time interval. The constructions are based on the detailed analysis of a low order Taylor expansion of the solution with a cleverly chosen initial datum. This has also been explored for the 2D Euler and 2D quasi-geostrophic equations by Li~\cite{li_euler_2012}, and for the 2D Boussinesq system by Zhang \cite{zhang_boussinesq_2013}, with similar results. 

Our objective in this article is to show the existence of solutions to the Navier--Stokes equations, both on $\mathbb{R}^2$ and on $\mathbb{T}^2$, which exhibit any (generic) bifurcation of stagnation points, generated by sequences of mergings and splittings of pairs of stagnation points. Apart from the obvious topological obstructions of such bifurcations, there are no restrictions on the number of points involved, on their positions, or on their type. Bifurcations of the kind constructed by Moffatt and Li--Sinai arise as particular cases. The strength of our approach is that it does not rely on the analysis of carefully constructed initial data, but rather on soft, general arguments.

To make this precise, it is convenient to introduce the {\em stream function}\/ $\psi(\bx,t)$, which determines the fluid velocity field through its perpendicular gradient:
\begin{equation}\label{eq:velocityfield}
    \bu=\nabla^\perp \psi:=\left(\frac{\partial \psi}{\partial y},-\frac{\partial\psi}{\partial x}\right).
\end{equation}
As is well known, the 2D Navier--Stokes equations can be written in terms of~$\psi$ as the scalar nonlocal  equation
\begin{equation}\label{eq:navierstokes2d}
    \frac{\partial \psi}{\partial t} + \Delta^{-1}\left(\nabla^\perp\psi\cdot\nabla \Delta\psi\right)=\Delta \psi,
\end{equation}
where $\Delta^{-1}$ is the inverse of the Laplacian and where $\nabla$ denotes the spatial gradient.
It is classical that this equation is globally (forward) wellposed for smooth enough initial data
\[
\psi(\bx,0)=\psi_0(\bx),
\]
related to the initial velocity of the fluid as $\bu(\cdot,0)=\nabla^\perp\psi_0$.

In terms of the stream function, the stagnation points of the fluid at time~$t$ correspond to the (spatial) {\em critical points}\/ of~$\psi$, that is, the points where $\nabla\psi(\cdot, t)$ vanishes. Throughout the paper, we will omit the word ``spatial''. Therefore,  studying the evolution of the stagnation points of the fluid amounts to keeping track of the critical points of~$\psi(\cdot,t)$, understood as a family of functions on~$\mathbb{R}^2$ parametrized by the time variable $t\in\mathbb{R}_+$. The key element here is that $\psi(\cdot,t)$ is obtained by solving Navier--Stokes with some initial datum~$\psi_0$.

A critical point of~$\psi(\cdot,t)$ is called {\em nondegenerate}\/ if the Hessian $\nabla^2\psi(\cdot,t)$ at that point has full rank. In particular, this condition implies that the critical point is {\em isolated}\/: in a small neighborhood, $\psi(\cdot,t)$ does not have any other critical points. There are two different types of {nondegenerate} critical points: local \textit{extrema}\/, corresponding to local {maxima} or {minima} of~$\psi(\cdot,t)$, and {\em saddles}\/. Since $\bu(\cdot,t)$ is the Hamiltonian vector field generated by~$\psi(\cdot,t)$, the {\em streamlines}\/ of the fluid at time~$t$ (which are the integral curves of~$\bu(\cdot, t)$) near an extremum are elliptic. The corresponding stagnation points are therefore \textit{centers}, or \textit{eddies}, of the flow, with the fluid {swirling} clockwise around the maxima and counterclockwise around the minima~\cite{Mo}.

\begin{figure}
        \def\svgwidth{8cm}
        %% Creator: Inkscape 1.4.2 (f4327f4, 2025-05-13), www.inkscape.org
%% Accompanies image file 'merging.pdf' (pdf, eps, ps)
\begingroup%
  \makeatletter%
  \providecommand\color[2][]{%
    \errmessage{(Inkscape) Color is used for the text in Inkscape, but the package 'color.sty' is not loaded}%
    \renewcommand\color[2][]{}%
  }%
  \providecommand\transparent[1]{%
    \errmessage{(Inkscape) Transparency is used (non-zero) for the text in Inkscape, but the package 'transparent.sty' is not loaded}%
    \renewcommand\transparent[1]{}%
  }%
  \providecommand\rotatebox[2]{#2}%
  \newcommand*\fsize{\dimexpr\f@size pt\relax}%
  \newcommand*\lineheight[1]{\fontsize{\fsize}{#1\fsize}\selectfont}%
  \ifx\svgwidth\undefined%
    \setlength{\unitlength}{515.3368612bp}%
    \ifx\svgscale\undefined%
      \relax%
    \else%
      \setlength{\unitlength}{\unitlength * \real{\svgscale}}%
    \fi%
  \else%
    \setlength{\unitlength}{\svgwidth}%
  \fi%
  \global\let\svgwidth\undefined%
  \global\let\svgscale\undefined%
  \makeatother%
  \begin{picture}(1,0.57138787)%
    \lineheight{1}%
    \setlength\tabcolsep{0pt}%
    \put(0,0){\includegraphics[width=\unitlength,page=1]{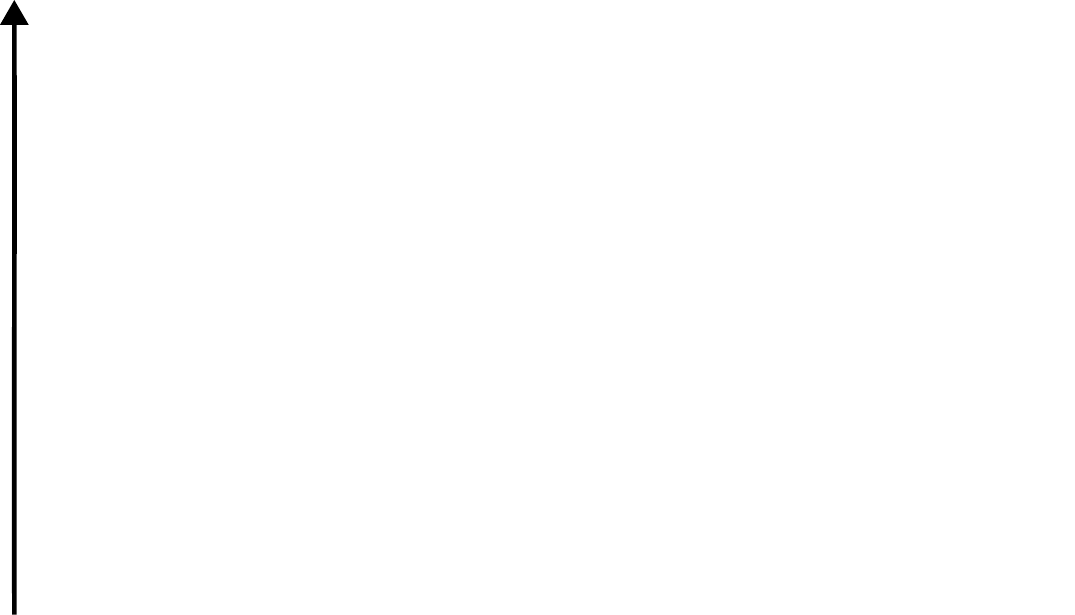}}%
    \put(0.03891878,0.53580862){\color[rgb]{0,0,0}\makebox(0,0)[lt]{\lineheight{1.25}\smash{\begin{tabular}[t]{l}$t$\end{tabular}}}}%
    \put(0.49934826,0.43654479){\color[rgb]{0,0,0}\makebox(0,0)[lt]{\lineheight{1.25}\smash{\begin{tabular}[t]{l}$\Omega$\end{tabular}}}}%
    \put(0.86351693,0.19675245){\color[rgb]{0,0,0}\makebox(0,0)[lt]{\lineheight{1.25}\smash{\begin{tabular}[t]{l}$t_1$\end{tabular}}}}%
    \put(0.86374189,0.33525508){\color[rgb]{0,0,0}\makebox(0,0)[lt]{\lineheight{1.25}\smash{\begin{tabular}[t]{l}$t_c$\end{tabular}}}}%
    \put(0.86362255,0.47870474){\color[rgb]{0,0,0}\makebox(0,0)[lt]{\lineheight{1.25}\smash{\begin{tabular}[t]{l}$t_2$\end{tabular}}}}%
    \put(0,0){\includegraphics[width=\unitlength,page=2]{merging.pdf}}%
  \end{picture}%
\endgroup%

        \caption{Merging of stagnation points: at time~$t_1$, the fluid exhibits two stagnation points (in red and blue), which merge into one (in green) at a later time~$t_c$, and subsequently disappear. Splitting corresponds to reflecting the curve across a horizontal plane.}
        \label{fig:merging}
\end{figure}

One says that a \textit{bifurcation} takes place at a time~$t$ if any critical points appear or disappear, or if any of those already existing changes its topological {type}. The simplest versatile scenario for the splitting and merging of stagnation points corresponds to the {\em saddle-node bifurcation}\/~\cite{Mo}, which involves the collapse of a center and a saddle. More precisely:
\begin{itemize}
    \item \textit{Merging} happens when a center and a saddle collide into a degenerate critical point at time $t$, which then disappears.
    \item \textit{Splitting} occurs when a new degenerate critical point  appears at time $t$, which then branches into a center and a saddle.
\end{itemize}
The splitting and merging of stagnation points is depicted schematically in Figure~\ref{fig:merging}. 

It is clear from Figure~\ref{fig:merging} that the splitting and merging of stagnation points should be analyzed in spacetime $\mathbb{R}^2\times\mathbb{R}_+$. Informally, what the figure shows is that stagnation points that eventually merge should be regarded as lying on a spacetime curve, or ``trajectory'', which looks like a cap; the merging occurs precisely at the top of the cap. Likewise, splitting stagnation points lie on trajectories that are shaped like a cup and the splitting happens at the bottom. 

To make use of this observation, let us introduce the following definition, where the {\em height}\/ function $h:\mathbb{R}^2\times\mathbb{R}_+\to\mathbb{R}$ is simply the projection on the time axis, $h(\bx,t):=t$. Basically, we are making sure that the points sitting on the two endpoints of a ``merging arc'' do merge at the point where the curve attains its maximum, as in the figure. 
\begin{definition*}
    A {\em merging (resp., splitting) arc}\/ is a smooth parametrized curve in spacetime, $\Gamma:[0,1]\to \mathbb{R}^2\times\mathbb{R}_+$, such that the function $h\circ\Gamma:[0,1]\to\mathbb{R}$ is strictly concave (resp., strictly convex) and attains its maximum (resp., minimum) in the interior $(0,1)$.

    By a {\em connecting arc} we refer to either a splitting or a merging arc. With some abuse of notation, in some parts of this article we will denote the image curve $\Gamma([0,1])$ also by $\Gamma$.
\end{definition*}

\begin{figure}
        \def\svgwidth{14cm}
        \hspace*{2.75cm}
        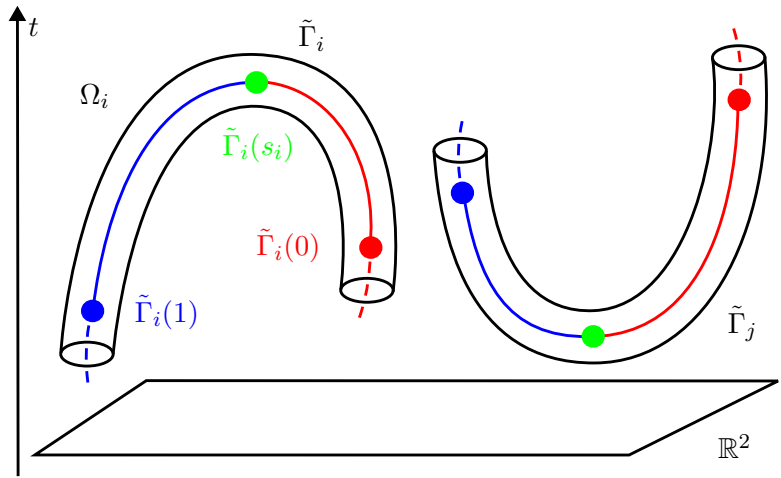
        \caption{Red stagnation points represent eddies (nodes) and blue ones represent saddles; on the left a saddle and a node merge into a degenerate stagnation point that immediately disappears, while on the right a degenerate stagnation point is created and then splits into a saddle and a node.}
        \label{fig:maintheoremobjects}
\end{figure}

We are now ready to state our main result. It ensures that, given any finite collection of splitting and merging arcs, as shown in Figure~\ref{fig:maintheoremobjects}, there exists a solution to the Navier--Stokes equations on $\mathbb{R}^2$ realizing that bifurcation pattern, up to a small spacetime diffeomorphism (which does not change the topology of the bifurcations, but which makes the trajectory that the fluid stagnation points follow slightly different). Note that the fluid may exhibit other stagnation points besides those described in the statement.

\begin{theorem}\label{thm:mergingsplittingR2}
    Let us fix a positive integer $r\geq 1$ and some small $\varepsilon>0$. Consider $k$ connecting arcs $\Gamma_i$, which we can assume to be merging for $1\leq i\leq k'$ and splitting for $k'<i\leq k$. 
    
    Then there exists some initial datum $\psi_0\in C^\infty_c(\mathbb{R}^2)$ for which the associated solution $\psi:\mathbb{R}^2\times\mathbb{R}_+\to\mathbb{R} $ to the Navier--Stokes equations realizes that bifurcation pattern, up to a small deformation. More precisely, there exists a diffeomorphism $\Phi$ of~$\mathbb{R}^2\times\mathbb{R}_+$ with 
$\|\Phi-\mathbf{id}\|_{C^r(\mathbb{R}^2\times\mathbb{R}_+)}<\varepsilon$ such that the slightly deformed spacetime curves $\widetilde{\Gamma}_i:=\Phi(\Gamma_i)$ consist of isolated stagnation points of the fluid. %\textcolor{blue}{That is, $\nabla \psi\circ \widetilde{\Gamma}_i\equiv 0$ and there exists an open tubular neighborhood $\Omega_i\subset\mathbb{R}^2\times\mathbb{R}_+$ whose subset of critical points of $\psi$ is the trace of $\tilde{\Gamma}_i$.}

Furthermore, the functions $h\circ\widetilde{\Gamma_i}:[0,1]\to\mathbb{R}$ are still strictly concave or convex, and attain their maximum or minimum at a unique interior point $s_i\in(0,1)$. The splitting (for $i>k'$) or merging (for $i\leq k'$) of stagnation points happens at the spacetime point $\widetilde{\Gamma}_i(s_i)$, which is a degenerate critical point of~$\psi$. The stagnation point $\widetilde{\Gamma}_i(s)$ can be chosen to be a nondegenerate local extremum for $s\in[0,s_i)$ and a nondegenerate saddle for $s\in(s_i,1]$.
\end{theorem}

\begin{remark}
    One can freely specify whether each extremum is a local maximum or minimum (so that they are clockwise or counterclockwise eddies), and one can also reverse the location of the extremum and the saddle along each curve $\widetilde{\Gamma}_i$. What is important is that each curve connects one saddle and one center.
\end{remark}

\begin{remark}
    The proof ensures that the bifurcation described in Theorem \ref{thm:mergingsplittingR2} is \textit{structurally stable}: any function that is close enough to $\psi$ in the $C^2$ norm, in a neighborhood of the connecting curves, presents the same splittings and mergings as those in Theorem~\ref{thm:mergingsplittingR2} up to a diffeomorphism.
\end{remark}
\begin{remark}
In Section \ref{section:torus} we prove an analogous result for the case of the torus $\mathbb{T}^2$. The only difference is that the analogue of the diffeomorphism $\Phi$ on $\mathbb{T}^2\times\mathbb{R}_+$ is no longer close to the identity (although one has quite a lot of information about its structure). 
\end{remark}

The proof of Theorem \ref{thm:mergingsplittingR2} is presented in Section \ref{section:proofmain}; in fact we prove a considerably more general statement concerning the realization of sets of stagnation points, see Theorem~\ref{thm:maintheoremR2}. All the technical details are relegated to Sections \ref{section:surface} and \ref{section:approx}, with some additional comments in Appendices \ref{app:appendixA} and \ref{app:appendixB}. Lastly, as mentioned before, in Section \ref{section:torus} we state an analogous version of the main theorem for the case of the torus $\mathbb{T}^2$, and we sketch its proof, which is similar to the case of $\mathbb{R}^2$.

\section{Proof of the main theorem: Links realized as sets of stagnation points}\label{section:proofmain}

Our main theorem is a consequence of the more general result that every compact curve in $\mathbb{R}^2\times\mathbb{R}_+$, up to a diffeomorphism close to the identity, can be realized as an isolated subset of critical points of a solution $\psi$ to the Navier-Stokes equations in $\mathbb{R}^2$. Each of its interior height critical points corresponds to a merging or splitting of stagnation points; the former in the case where it is a maximum and the latter if it is a minimum. If the curve is \textit{closed} (a knot), i.e., if it does not have endpoints, then the total number of critical points of $\psi$ that split coincides with the number of the ones that merge. This in turn can be generalized to any finite family of compact curves in $\mathbb{R}^2\times\mathbb{R}_+$, some of them without endpoints and others not.

\Huge

\begin{figure}
        \def\svgwidth{10cm}
        \hspace*{1cm}
        %% Creator: Inkscape 1.4.2 (f4327f4, 2025-05-13), www.inkscape.org
%% Accompanies image file 'openlink.pdf' (pdf, eps, ps)
\begingroup%
  \makeatletter%
  \providecommand\color[2][]{%
    \errmessage{(Inkscape) Color is used for the text in Inkscape, but the package 'color.sty' is not loaded}%
    \renewcommand\color[2][]{}%
  }%
  \providecommand\transparent[1]{%
    \errmessage{(Inkscape) Transparency is used (non-zero) for the text in Inkscape, but the package 'transparent.sty' is not loaded}%
    \renewcommand\transparent[1]{}%
  }%
  \providecommand\rotatebox[2]{#2}%
  \newcommand*\fsize{\dimexpr\f@size pt\relax}%
  \newcommand*\lineheight[1]{\fontsize{\fsize}{#1\fsize}\selectfont}%
  \ifx\svgwidth\undefined%
    \setlength{\unitlength}{525.94199834bp}%
    \ifx\svgscale\undefined%
      \relax%
    \else%
      \setlength{\unitlength}{\unitlength * \real{\svgscale}}%
    \fi%
  \else%
    \setlength{\unitlength}{\svgwidth}%
  \fi%
  \global\let\svgwidth\undefined%
  \global\let\svgscale\undefined%
  \makeatother%
  \begin{picture}(1,0.56178081)%
    \lineheight{1}%
    \setlength\tabcolsep{0pt}%
    \put(0,0){\includegraphics[width=\unitlength,page=1]{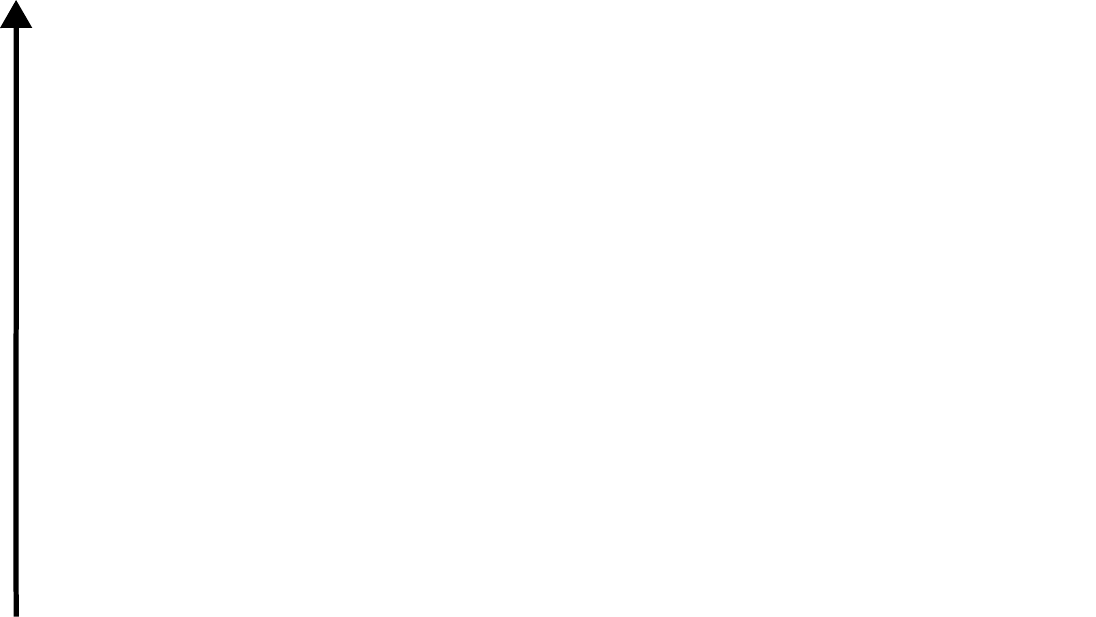}}%
    \put(0.03979082,0.52500508){\color[rgb]{0,0,0}\makebox(0,0)[lt]{\lineheight{1.25}\smash{\begin{tabular}[t]{l}$t$\end{tabular}}}}%
    \put(0.81948534,0.24826629){\color[rgb]{0,0,0}\makebox(0,0)[lt]{\lineheight{1.25}\smash{\begin{tabular}[t]{l}$\zeta$\end{tabular}}}}%
    \put(0,0){\includegraphics[width=\unitlength,page=2]{openlink.pdf}}%
    \put(0.33115462,0.48286758){\color[rgb]{0,0,1}\makebox(0,0)[lt]{\lineheight{1.25}\smash{\begin{tabular}[t]{l}$L^-$\end{tabular}}}}%
    \put(0.67861089,0.3783582){\color[rgb]{1,0,0}\makebox(0,0)[lt]{\lineheight{1.25}\smash{\begin{tabular}[t]{l}$L^+$\end{tabular}}}}%
  \end{picture}%
\endgroup%

        \caption{A pair of compact curves realized by stagnation points. The arrows indicate the orientation of the curves.}
        \label{fig:openlink}
\end{figure}

\normalsize

 Let $L$ be the union of a finite family of {$m$ disjoint curves in $\mathbb{R}^2\times\mathbb{R}_+$, each diffeomorphic to either $[0,1]$ or $\mathbb{T}^1:=\mathbb{R}/2\pi\mathbb{Z}$. If all the components of $L$ are diffeomorphic to the latter, then $L$ is a \textit{link}. Otherwise, it is an \textit{open link}, as in Figure \ref{fig:openlink}.} The curve depicted on the right side is a trefoil knot, and the one to the left is diffeomorphic to $[0,1]$. Suppose that $L$ is in general position with respect to the height function, i.e., the critical points of the restriction $h|L$ are all nondegenerate, and fix an orientation $\zeta$ in $L$.

The \textit{horizontal points} of $L$ are the critical points of $h|L$.  For simplicity, if $L$ is an open link, assume that none of its endpoints are horizontal points. A \textit{positive interval} of $L$ is an open segment between two of its horizontal points, or between and endpoint and a horizontal point in the case of an open link, along which $h$ increases following the orientation prescribed by $\zeta$. Likewise, a \textit{negative interval} of $L$ is one for which $h$ decreases following the aforementioned orientation. The union of the former type of intervals is denoted as $L^+$, while we use $L^-$ for the latter. Positive and negative intervals are swapped under the reversal $\zeta\mapsto-\zeta$. Both, together with the horizontal points (and the endpoints in the case of an open link), form a {partition} of $L$. In Figure \ref{fig:openlink}, the positive intervals are depicted in red, the negative ones in blue, and the horizontal points in green. The orientation $\zeta$ is given by the different arrows. Observe that, since $L$ is in general position with respect to the height function, there is only a finite quantity of these three types of components and the horizontal points coincide with the height maxima or minima of the curve. Now, let us choose a finite subcollection of positive intervals, denote its union as $L^+_{\max}$, and define the complement $L^+_{\min}:=L^+\backslash L^+_{\max}$.

\begin{theorem}\label{thm:maintheoremR2}
    Fix a positive integer $r\geq 1$ and let $\varepsilon>0$ be any real number. Let $L$ be an arbitrary link in $\mathbb{R}^2\times\mathbb{R}_
    +$, possibly open, with an orientation $\zeta$. There exists a solution $\psi:\mathbb{R}^2\times\mathbb{R}_+\rightarrow\mathbb{R}$ to the 2D Navier-Stokes equations, a tubular neighborhood $\Omega$ of $L$, and a diffeomorphism $\Phi:\Omega\rightarrow\Omega$ such that:
    \begin{itemize}
        \item $\|\Phi-\mathbf{id}\|_{C^r(\Omega)}<\varepsilon$.
        \item $L_2:=\Phi(L)$ is an isolated subset of the set of critical points of $\psi$ in $\mathbb{R}^2\times\mathbb{R}_+$. That is, the set of critical points of $\psi$ in $\Omega$ is exactly $L_2$.
    \end{itemize}
    Moreover, the type of each critical point of $\psi$ along $L_2$ is known.
    \begin{itemize}
        \item Consider the orientation induced by $\zeta$ in $L_2$ via $\Phi$. Then, the critical points of $\psi$ in the positive intervals $L_2^+$ are centers, i.e., extrema, and those in $L_2^-$ are saddles.
        \item Consider the decomposition of $L_2^+$ into $L_{2,\max}^+$ and $L_{2,\min}^+$ coming from $L^+$. The critical points of $\psi$ in the former intervals are maxima, and those in the latter are minima.
    \end{itemize}
\end{theorem}
\noindent As in Theorem \ref{thm:mergingsplittingR2}, the situation depicted by Theorem \ref{thm:maintheoremR2} is structurally stable. Any stream function that is close enough to $\psi$ in the $C^r$ norm also realizes $L$ as a subset of its critical points up to a diffeomorphism. At each height maximum of $L_2$ there is a merging of critical points, and at each minimum a splitting. So, it is clear that Theorem \ref{thm:maintheoremR2} extends the geometrical picture depicted in Theorem \ref{thm:mergingsplittingR2}, which follows from it as a straightforward consequence. Now, several remarks are in order.

\begin{remarks}\label{remarks:commentstheoremR2}
    First, the regularity of $L$ in Theorem \ref{thm:maintheoremR2} and of the connecting arcs in Theorem \ref{thm:mergingsplittingR2} can be weakened: we only need $C^k$ with $k\geq 0$. Indeed, the first step of the proof of Theorem \ref{thm:maintheoremR2} consists in approximating $L$ by an analytic link, which can be done for any regularity class thanks to Whitney's Analytic Approximation Theorem \cite[Theorem 5.6.]{hirsch_differential_1976}.
    
    Second, let us explain in detail how the partition of $L$ descends to one of $L_2$ in the latter half of Theorem \ref{thm:maintheoremR2}. The orientation considered in $L_2$ is, between any two of its points, the one prescribed by $\zeta$ between their two inverse images in $L$ via $\Phi^{-1}$. For the partition of $L_2^+$ into $L_{2,\max}^+$ and $L_{2,\min}^+$, we consider that $\varepsilon$ in the theorem's statement is small enough so that each positive interval in $L_2^+$ intersects just one of the components of $\Phi(L^+)$. We define $L_{2,\max}^+$ as those that overlap with $\Phi(L_{\max}^+)$, and choose $L_{2,\min}^+$ as those positive intervals whose intersection with $\Phi(L_{\min}^+)$ is nonempty.

    Third, we insist that the result of Theorem \ref{thm:maintheoremR2}, and thus of Theorem \ref{thm:mergingsplittingR2}, is contingent on deforming $L$ via $\Phi$. Our procedure does not allow to directly realize $L$ as a subset of the set of critical points of a stream function $\psi$ without passing through a diffeomorphism first. This deformation step, in fact, may even be a necessary condition in some cases, see e.g. \cite[Remark 3.13]{enciso_topological_2018} about how it is essential that plane curves are defined modulo diffeomorphism to realize the level sets of harmonic functions in $\mathbb{R}^2$.

    Fourth, technically the link $L$ is realized as an isolated subset of critical points of $\psi$ if and only if it is closed, i.e., if it is not an open link. Otherwise, the set of critical points of $\psi$ extend the interval-like components of $L$ into closed curves {diffeomorphic to $\mathbb{T}^1$}. This is an artifact coming from the proof of the theorem itself, whose first step consists in extending $L$ so it is a (closed) link. The theorem's statement holds because afterwards we can reduce $\Omega$ so it does not contain any points of the new extension, ``cutting'' where the endpoints of $L$ are.
\end{remarks}

An analogous result to Theorem \ref{thm:maintheoremR2} exists for the case of the torus $\mathbb{T}^2$. As with Theorem \ref{thm:mergingsplittingR2}, it is deferred to Section \ref{section:torus}. The only difference is that, in $\mathbb{T}^2$, we have to scale and {compress} the link $L$ before it is realized as a subset of stagnation points.

\subsection{Proof of Theorem~\ref{thm:maintheoremR2}}
Since Theorem \ref{thm:mergingsplittingR2} is a consequence of Theorem~\ref{thm:maintheoremR2}, we first prove the latter. The strategy consists of three steps:
\begin{enumerate}
    \item First, we study the mergings and splittings of stagnation points for the linearization of the two-dimensional Navier-Stokes equations at the trivial solution, which turns out to be the heat equation. We then prove an analog version of Theorem \ref{thm:maintheoremR2} for this case, which makes use of the Cauchy-Kovalevskaya theorem and a global approximation result for parabolic PDEs. The resulting heat solutions are robust against perturbations.
    \item Second, using a suitable scaling, we construct solutions to Navier-Stokes \eqref{eq:navierstokes2d} in $\mathbb{R}^2\times\mathbb{R}_+$ which are $C^r$-close to the heat solution previously constructed in the domain $\Omega$. This is done via an elementary perturbative-type argument.
    \item Finally, we transfer our results about mergings and splittings for the robust solutions to the heat equation found in the first step of the proof to solutions of the Navier-Stokes equations in $\mathbb{R}^2$.
\end{enumerate}
Let $L$ be a link as in the statement of Theorem~\ref{thm:maintheoremR2}, partitioned into positive and negative intervals through an orientation $\zeta$. We can safely assume that all its components are knots, i.e., that it is not an open link. Otherwise, extend it so it has no endpoints, all while staying in the $C^\infty$ class and in the domain $\mathbb{R}^2\times\mathbb{R}_+$. Fix an integer $r\geq 1$ and a real number $\varepsilon>0$.

Obviously, the linearization of Equation~\eqref{eq:navierstokes2d} around its trivial solution is the heat equation, which for $v=v(x,y,t)$ we write as:
\begin{equation}\label{eq:heatequation}
    \frac{\partial v}{\partial t}=\Delta v,
\end{equation}
where $\Delta:=\partial_{xx}+\partial_{yy}$ denotes the spatial Laplacian. As before, the \textit{critical points} for $v$ are the points in $\mathbb{R}^2\times\mathbb{R}_+$ where the two partial derivatives $\partial_x v$ and $\partial_y v$ are zero.

We want to show that, up to a diffeomorphism $C^r$-close to the identity, the link $L$ is realized as an isolated subset of critical points of a solution to \eqref{eq:heatequation}. This has to be done in a way that is robust under pertubations, so that $L$ is preserved when passing from the linear equation \eqref{eq:heatequation} back to \eqref{eq:navierstokes2d}. More precisely, we need to construct a solution $v:\mathbb R^2\times \mathbb R_+\to \mathbb R$ to \eqref{eq:heatequation} so that:
\begin{itemize}
    \item[\textbf{(A)}] \textit{The link $L$ is realized as a subset of critical points of $v$.} The two derivatives $\partial_x v$ and $\partial_y v$ have to be identically zero on the image of $L$ via a diffeomorphism $\Phi_1:\Omega\rightarrow\Omega$, where $\Phi_1$ can be chosen arbitrarily close to the identity in the $C^r$ norm and $\Omega$ is a tubular neighborhood of $L$.

    \item[\textbf{(B)}] \textit{The type of the critical points of $v$ along $L$ is known.} After deforming $L$ via $\Phi_1$ and defining $L_1:=\Phi(L)$, if we consider the partition of $L_1$ into positive and negative intervals, the critical points of $v$ in $L^+_{1,\max}$ have to be maxima, the ones in $L^+_{1,\min}$ minima, and the ones in $L_1^-$ saddles.

    \item[\textbf{(C)}] \textit{The solution $v$ is structurally stable around $L$.} The zero set of the gradient field $(v_x,v_y)$, where $v_x:=\partial_x v$ and $v_y:=\partial_y v$, has to be robust in a neighborhood of $L_1$, which in our case is tantamount to demanding that:
\begin{equation}\label{eq:structuralstability}
        \rank\big(\nabla_3 v_x(\xi),\nabla_3 v_y(\xi)\big)=2,\quad \forall \xi:=(x,y,t) \in L_1.
\end{equation}
Here and in what follows, $\nabla_3$ is the gradient in the (x,y,t)-variables, and so $\big(\nabla_3 v_x(\xi),\nabla_3 v_y(\xi)\big)$ is interpreted as a $2\times 3$ matrix.
\end{itemize}

The purpose of Conditions \textbf{(A)} and \textbf{(B)} is clear. Condition \textbf{(C)} is what guarantees that the two properties \textbf{(A)} and \textbf{(B)} are preserved when passing from $v$ to $\psi$ in the next step of the proof. The existence of a solution $v$ to the heat equation satisfying these three conditions is established in the following:

\begin{proposition}\label{prop:heat}
    Let $\varepsilon_1>0$ be any real number. There exists a tubular neighborhood $\Omega$ of $L$, a diffeomorphism $\Phi_1:\Omega\rightarrow\Omega$, and a local solution to the heat equation $v:\Omega\rightarrow\mathbb{R}$ such that:
    \begin{enumerate}
        \item $\|\Phi_1-\mathbf{id}\|_{C^r(\Omega)}<\varepsilon_1$.
        \item The local solution $v$ satisfies Conditions \textbf{(A)}, \textbf{(B)} and \textbf{(C)}.
    \end{enumerate}
\end{proposition}

The proof of this proposition is presented in Section \ref{section:surface}. It consists in constructing a Cauchy problem for the heat equation such that, by adequately choosing the Cauchy surface and data, any solution is forced to satisfy \textbf{(A)}, \textbf{(B)} and \textbf{(C)}. If both the surface and the data are analytic, then such a solution exists by the {Cauchy-Kovalevskaya Theorem} \cite[Proposition 6.4.2.]{taylor_partial_2011}. It is important to note then that the diffeomorphism $\Phi_1$ in Proposition \ref{prop:heat} only serves to change $L$ into a $C^\omega$ curve. If $L$ is already analytic beforehand, then $\Phi_1$ can be chosen as the identity.

The solution $v$ in Proposition~\ref{prop:heat} serves as a first approximation to the Navier-Stokes stream function $\psi$ in Theorem \ref{thm:maintheoremR2}. There are two difficulties to overcome:
\begin{itemize}
    \item First, the function $v$ is only defined locally, while we want $\psi$ to be defined in all $\mathbb{R}^2\times\mathbb{R}_+$. A way to solve this is to pass from $v$ to a global heat solution, which we can do with one of the {approximation theorems for parabolic equations} in \cite{enciso_approximation_2019}.
    \item Second, we need a way to control the nonlinear term in \eqref{eq:navierstokes2d} to promote the previous global heat solution to a Navier-Stokes stream function. This can be done using a {scaled two-dimensional Navier-Stokes equations} for $\tilde{\psi}=\tilde{\psi}(x,y,t)$ as:
\begin{equation}\label{eq:scalednavierstokes2d}
    \frac{\partial \tilde{\psi}}{\partial t} + \delta\Delta^{-1}\left(\frac{\partial\tilde{\psi}}{\partial x}\cdot\frac{\partial (\Delta\tilde{\psi})}{\partial y}-\frac{\partial(\Delta\tilde{\psi})}{\partial x}\cdot\frac{\partial\tilde{\psi}}{\partial y}\right)=\Delta \tilde{\psi},
\end{equation}
where $\delta>0$ is the {scaling parameter}. It is clear that $\psi(x,y,t):=\delta\tilde{\psi}(x,y,t)$ is a solution to the (nonscaled) Navier-Stokes equations \eqref{eq:navierstokes2d}.
\end{itemize}

Combining these ideas, we find solutions $\tilde{\psi}$ to Equation~\eqref{eq:scalednavierstokes2d} for small $\delta$ which approximate the local heat solution $v$:

\begin{proposition}\label{prop:approx}
    Let $\delta>0$ be any real number, and suppose that $v:\Omega\rightarrow\mathbb{R}$ is the solution to the heat equation in Proposition \ref{prop:heat}. Reducing $\Omega$ if necessary, there exists a solution $\tilde{\psi}:\mathbb{R}^2\times\mathbb{R}_+\rightarrow\mathbb{R}$ to Equation~\eqref{eq:scalednavierstokes2d} in $\mathbb{R}^2$, with scaling parameter $\delta$, such that $\|v-\tilde{\psi}\|_{C^{r+1}(\Omega)}<\delta$.
\end{proposition}

The proof of this proposition is presented in Section \ref{section:approx}. It follows the two-step procedure previously outlined, and is similar to the one used in \cite[Section 5]{enciso_approximation_2021} for the case of the Schrödinger equation.

To complete the proof of the theorem, take the solution $\tilde{\psi}:\mathbb{R}^2\times\mathbb{R}_+\rightarrow\mathbb{R}$ from Proposition \ref{prop:approx} and let $\varepsilon_2>0$ be any real number. The $C^{r+1}$ norm of the difference between $v$ and $\tilde{\psi}$ is bounded by $\delta$, so:
\begin{equation}
    \Big\|\Big(\frac{\partial\tilde{\psi}}{\partial x},\frac{\partial\tilde{\psi}}{\partial y}\Big)-\Big(\frac{\partial v}{\partial x},\frac{\partial v}{\partial y}\Big)\Big\|_{C^r(\Omega)}<\delta.
\end{equation}
Since $v$ satisfies Condition \textbf{(C)}, we can choose $\delta$ small enough so that $\tilde{\psi}$ inherits \textbf{(A)} and \textbf{(B)}. Indeed, by {Thom's Isotopy Theorem} \cite[Theorem 3.14.]{enciso_topological_2018}, there exists a diffeomorphism $\Phi_2:\Omega\rightarrow\Omega$ which is $\varepsilon_2$-close to the identity in the $C^r$ norm such that the set of critical points of $\tilde{\psi}$ in $\Omega$ is $L_2:=\Phi_2(L_1)=(\Phi_2\circ\Phi_1)(L)$. Moreover, the types of the critical points that $\tilde{\psi}$ has along $L_2$ coincides with those that $v$ has along $L_1$, which thanks to Condition \textbf{(B)} is in correspondence with those chosen for $L$ in the statement of Theorem \ref{thm:maintheoremR2}.

To finish the proof, first consider the diffeomorphism $\Phi:=\Phi_2\circ\Phi_1:\Omega\rightarrow\Omega$. By reducing $\varepsilon_1$ and $\varepsilon_2$, it satisfies $\|\Phi-\mathbf{id}\|<\varepsilon$. Second, let us define $\psi:=\delta\tilde{\psi}:\mathbb{R}^2\times\mathbb{R}_+\rightarrow\mathbb{R}$. It is a solution to the  (nonscaled) Navier-Stokes equations in $\mathbb{R}^2$. Its critical points are given by $\partial_x \psi=\partial_y \psi=0$, and since $\psi$ is proportional to $\tilde\psi$, the sets of critical points of both functions coincide. Since $\delta>0$, the type of said points is the same too, and Theorem \ref{thm:maintheoremR2} then follows.

\subsection{Proof of Theorem \ref{thm:mergingsplittingR2}.} As stated before, Theorem \ref{thm:mergingsplittingR2} can be obtained as a direct consequence of Theorem \ref{thm:maintheoremR2}. We just have to consider the open link $L$ made of the connecting arcs $\Gamma_i$. To match the choice of saddles and extrema made in Theorem \ref{thm:mergingsplittingR2}, we choose the orientation $\zeta$ that {goes from $\Gamma_i(0)$ to $\Gamma_i(1)$ for $1\leq i\leq k^\prime$, and from $\Gamma_i(1)$ to $\Gamma_i(0)$ for $k^\prime<i\leq k$}. The result is then straightforward.

\section{Proof of Proposition~\ref{prop:heat} }\label{section:surface}
 
In this section we present the proof of Proposition \ref{prop:heat}. The main idea is to construct a Cauchy problem for the heat equation tailored in a such way that all of its solutions $v$ have to satisfy the three conditions \textbf{(A)}, \textbf{(B)} and \textbf{(C)}.

Let $L$ be a link partitioned into positive and negative intervals through an orientation $\zeta$ as in Section \ref{section:proofmain}, and suppose that $\Sigma\subset\mathbb{R}^3$ is a smooth orientable surface such that $L\subset \Sigma$. Let $\mathbf{N}:\Sigma\rightarrow\mathbb{S}^2$ be a {Gauss map} for $\Sigma$, {where $\mathbb{S}^2$ is the unit 2-sphere. Equivalently, $\mathbf{N}$ is a unitary normal field on the surface $\Sigma$.} Consider the Cauchy problem:
\begin{equation}\label{eq:heatcauchy}
    \begin{cases}
        \partial_t v=\Delta v,\\
        v|\Sigma=f,\,\mathbf{N}\cdot\nabla_3 v=g.
    \end{cases}
\end{equation}
The two functions $f,g:\Omega\rightarrow\mathbb{R}$ are the Cauchy data for the problem, with $\Omega$ a neighborhood of~$\Sigma$ (although, of course, the Cauchy problem only uses the values of $f$ and $g$ on~$\Sigma$). After possibly deforming $L$, if both the surface $\Sigma$ and the data are analytic, and $\Sigma$ is {noncharacteristic} for the heat equation, then by the {Cauchy-Kovalevskaya Theorem} \cite[Proposition 6.4.2.]{taylor_partial_2011} there exists a unique local solution $v:\Omega\rightarrow\mathbb{R}$, where $\Omega$ here is a subset of the original domain of the data.

With this in mind, our proof of Proposition \ref{prop:heat} consists in showing that, up to a diffeomorphism of $L$ arbitrarily close to the identity in the $C^r$ norm, it is possible to choose the surface $\Sigma$ and the data $f,g$ in a way that forces the local solution $v$ to satisfy \textbf{(A)}, \textbf{(B)} and \textbf{(C)}. The construction we present is completely explicit, starting first with the surface $\Sigma$, and then going on to the boundary data $f$ and $g$.

\subsection{Part I: Construction of a Noncharacteristic Boundary Surface}\label{subs:constsigma}

Since Proposition \ref{prop:heat} is local in nature, we can suppose that the link $L$ has only {one component $L_0$} (a smooth knot in $\mathbb{R}^2\times\mathbb{R}_+$) in general position with respect to the height function $h(x,y,t)=t$. In what follows, we use $\xi:=(x,y,t)$ to denote points in $\mathbb{R}^2\times\mathbb{R}_+$.

By {Whitney's Analytic Approximation Theorem} \cite[Theorem 5.6.]{hirsch_differential_1976} there exists a diffeomorphism $\Phi_1:\mathbb{R}^3\rightarrow\mathbb{R}^3$ with $\|\Phi_1-\mathbf{id}\|_{C^r(\mathbb{R}^3)}<\varepsilon_1$ such that {the image $L_1:=\Phi_1(L_0)$} is an analytic curve.

Let $\phi:\mathbb{S}^1\rightarrow\mathbb R^3$ be the arc-length parametrization of $L_1$, where $\mathbb{S}^1:=\mathbb{R}/(l\mathbb{Z})$ for $l$ the length of $L_1$. For any $\theta\in\mathbb{S}^1$
we take the unit {tangent vector} $\mathbf{t}(\theta)\in\mathbb{R}^3$ to $L_1$ at the point $\phi(\theta)=\xi$. We can assume that the orientation $\zeta$ in $L_1$ coincides with the orientation of the parametrization $\phi$. Then, let $L_1^+$ and $L_1^-$ be the corresponding positive and negative intervals of $L_1$, and make a choice of $L_{1,\max}^+$ and $L_{1,\min}^+$ as in Theorem \ref{thm:maintheoremR2}.
   
Consider the horizontal points of $L_1$, i.e., those $\xi=\phi(\theta)\in L_1$ for which $\mathbf{t}(\theta)\cdot\mathbf{e_3}=0$. By general position, there is only a finite quantity of them. Assume we denote their number as $m\in\mathbb{Z}$, and observe that $m\geq 2$ (there is at least one height minimum and one maximum). Let us denote them as $\xi_j\in L_1$ with $j\in \mathbb{Z}/ m\mathbb{Z}$, the labeling made according to the orientation prescribed by $\phi$ on $L_1$.
For all $j\in\mathbb{Z}/m\mathbb{Z}$, let $\theta_j\in\mathbb{S}^1$ be such that $\phi(\theta_j)=\xi_j$.

Let us construct an analytic surface $\Sigma\supset L_1$ that is noncharacteristic for the Cauchy problem \eqref{eq:heatcauchy}. To this end, let $\mathbf{k}:\mathbb{S}^1\rightarrow\mathbb{R}^3$ be an analytic map such that $\mathbf{k}(\theta)\cdot\mathbf{t}(\theta)=0$ and $\mathbf{k}(\theta)\cdot\mathbf{k}(\theta)=1$ for all $\theta\in\mathbb{S}^1$. In other words, $\mathbf{k}$ is a $C^\omega$ unit normal field in $L_1$. The surface $\Sigma$ will be obtained by ``sliding'' the points of $L_1$ along the direction of $\mathbf{k}$, which we appropriately call a \textit{sliding vector}. More precisely, consider the map:
\begin{align}\label{eq:sigmaparam}
    \begin{split}
        \Theta_0:\,&\mathbb{S}^1\times[-\sigma_1,\sigma_1]\rightarrow\mathbb{R}^3,\\
        &(\theta,\sigma)\mapsto \phi(\theta)+\sigma\mathbf{k}(\theta),
    \end{split}
\end{align}
where $\sigma_1>0$ is small enough. By the {Tubular Neighborhood Theorem} \cite[Theorem 6.4.]{lee_introduction_2012}, the map $\Theta_0$ is a diffeomorphism onto its image. We choose $\Sigma:=\Theta_0(\mathbb{S}^1\times[-\sigma_1,\sigma_1])$, which is analytic and diffeomorphic to a cylinder by construction and contains $L_1$. Next, we argue how to choose $\mathbf{k}$ so that it is noncharacteristic for the problem \eqref{eq:heatcauchy}.

At each point $\xi=\Theta_0(\theta,\sigma)\in\Sigma$, the two vectors 
\begin{align*}
&\mathbf{v}_{\theta,\Sigma}(\theta,\sigma):=\mathbf{t}(\theta)+\sigma\partial_\theta\mathbf{k}(\theta),\\
&\mathbf{v}_{\sigma,\Sigma}(\theta,\sigma):=\mathbf{k}(\theta),
\end{align*}
form a basis for the tangent space $T_{\xi}\Sigma$ as a subspace of $\mathbb{R}^3$. In addition, define:
\begin{align}\label{eq:normalsigma}
       \mathbf N(\theta,\sigma):= \frac{\mathbf{v}_{\theta,\Sigma}(\theta,\sigma)\times \mathbf{v}_{\sigma,\Sigma}(\theta,\sigma)}{\big\|\mathbf{v}_{\theta,\Sigma}(\theta,\sigma)\times \mathbf{v}_{\sigma,\Sigma}(\theta,\sigma)\big\|}.
\end{align}

It is easy to see that the surface $\Sigma$ is noncharacteristic for the Cauchy problem \eqref{eq:heatcauchy} if and only if:
\begin{equation}\label{eq:conds}
    \mathbf{N}(\theta,\sigma)\cdot\mathbf{e_3}\neq 0,\quad\forall(\theta,\sigma)\in\mathbb{S}^1\times[-\sigma_1,\sigma_1].
\end{equation}
By taking a smaller value of $\sigma_1$ if needed, the compactness of $L_1$ implies that we only need to make sure that \eqref{eq:conds} is satisfied on $L_1$, i.e.,  $\mathbf{N}(\theta,0)\cdot\mathbf{e_3}\neq 0$ for all $\theta\in\mathbb{S}^1$. This is automatically satisfied for those $\xi=\phi(\theta)\in L_1$ different from the horizontal points $\xi_j$. Consequently, by recalling \eqref{eq:normalsigma}, it all amounts to choosing $\mathbf{k}$ so that:
\begin{equation}\label{eq:condk}
    \mathbf{k}(\theta_j)\cdot\mathbf{e_3}\neq0,\quad\forall\,j=0,1,2\ldots,m-1\,\big(\text{mod}\,m\big).
\end{equation}
Of course, an analytic normal field $\mathbf{k}$ with this property exists, thus yielding the desired noncharacteristic surface $\Sigma$.

\subsection{Part II: Choice of Adapted Cauchy Data}

First we construct suitable coordinates parametrizing a tubular neighborhood $\Omega$ of $L_1$. Let $\rho_1>0$ be any real number, and consider the map:
\begin{align}\label{eq:omegaparam}
    \begin{split}
        \Theta_1:\,&\mathbb{S}^1\times[-\sigma_1,\sigma_1]\times[-\rho_1,\rho_1]\rightarrow\mathbb{R}^3,\\
        &(\theta,\sigma,\rho)\mapsto \phi(\theta)+\sigma\mathbf{k}(\theta)+\rho\mathbf{N}(\theta,\sigma).
    \end{split}
\end{align}
As before, if $\rho_1>0$ is small enough, $\Theta_1$ is a diffeomorphism onto its image. Set $W:=\mathbb{S}^1\times[-\sigma_1,\sigma_1]\times[-\rho_1,\rho_1]$ and define $\Omega:=\text{int }\Theta_1(W)$. By reducing $\sigma_1$ and $\rho_1$, we may suppose that $\overline{\Omega}\subset\mathbb{R}^2\times\mathbb{R}_+$. By construction, the map $\Theta_1$ is analytic, in particular, $\theta$, $\sigma$ and $\rho$ can be understood as analytic functions $\Omega\rightarrow\mathbb{R}$, and then, every map in $\Omega$ can be expressed in terms of them.

 At each point $\xi=\Theta_1(\theta,\sigma,\rho)\in\Omega$  we can write a basis for the tangent space of $\Omega$ at $\xi$:
\begin{equation}\label{eq:omegatangentbasis}
    \begin{cases}
        \mathbf{v}_{\theta,\Omega}(\theta,\sigma,\rho)=\mathbf{t}(\theta)+\sigma\partial_\theta\mathbf{k}(\theta)+\rho\partial_\theta\mathbf{N}(\theta,\sigma),\\
        \mathbf{v}_{\sigma,\Omega}(\theta,\sigma,\rho)=\mathbf{k}(\theta)+\rho\partial_\sigma\mathbf{N}(\theta,\sigma),\\
        \mathbf{v}_{\rho,\Omega}(\theta,\sigma,\rho)=\mathbf{N}(\theta,\sigma).
    \end{cases}
\end{equation}
Setting $\mathbf{n}(\theta):=\mathbf N(\theta,0)$, the set of triples
\begin{equation}\label{eq:movingframe}
    \mathcal{B}:=\big\{(\mathbf{t}(\theta),\mathbf{k}(\theta),\mathbf{n}(\theta)):\theta\in\mathbb{S}^1\big\},
\end{equation}
forms an orthonormal, positively oriented basis of $\mathbb{R}^3$ for each $\theta$, which is a {moving frame} for $L_1$.

We want to find analytic Cauchy data $f,g:\Omega\rightarrow\mathbb{R}$ such that any solution to \eqref{eq:heatcauchy} satisfies \textbf{(A)}, \textbf{(B)} and \textbf{(C)}. We divide this into three steps.
\begin{enumerate}
    \item First, we determine what requirements  Condition \textbf{(A)} imposes on $f$ and $g$. 
    \item Second, we present an {ansatz} for $f$ and $g$ satisfying the conditions found in the first step. 
    \item Third, we finish adjusting the data in order to satisfy the two remaining conditions \textbf{(B)} and \textbf{(C)}.
\end{enumerate}

\subsubsection{Step~1: System of equations imposed by the gradient nullity condition (A)}\label{sub:conditionAexplanation}

By taking a smaller $\Omega$ if needed, let $v:\Omega\rightarrow\mathbb{R}$ be a solution to \eqref{eq:heatcauchy}. 
We call $(x,y,t)$ the \textit{positional variables} and $(\theta,\sigma,\rho)$ the \textit{coordinate variables}. Both variables are related through the parametrization $\Theta_1:W\rightarrow\overline{\Omega}$.

By composing $v$ with $\Theta_1$, we can write $v$ as a function of $\theta$, $\sigma$ and $\rho$. Obviously, the partial derivatives of $v\circ\Theta_1$ with respect to $\theta$, $\sigma$ and $\rho$ correspond to the directional derivatives of $v$ along $\mathbf{v}_{\theta,\Omega}$, $\mathbf{v}_{\sigma,\Omega}$ and $\mathbf{v}_{\rho,\Omega}$, respectively. This way, we have
\begin{equation}\label{eq:dirderivatives}
    \partial_\theta v(\xi)=\nabla_3 v(\xi)\cdot\mathbf{v}_{\theta,\Omega}(\Theta_1^{-1}(\xi)),\quad\xi\in\Omega,
\end{equation}
and the same for the other two partial derivatives $\partial_\sigma$ and $\partial_\rho$. This definition can be generalized to higher-order derivatives, or even to partial derivatives in both the positional and the coordinate variables. This notation is explained using $v$ as an example, but it is valid for any function whose domain is $\Omega$ or a subset of it.

The three vectors $\mathbf{v}_{\theta,\Omega}$, $\mathbf{v}_{\sigma,\Omega}$ and $\mathbf{v}_{\rho,\Omega}$ depend solely on $\mathbf{t}$, $\mathbf{k}$ and their derivatives, as it can be seen in \eqref{eq:omegatangentbasis}. They are the duals of the gradients of the coordinates $\theta$, $\sigma$ and $\rho$ when seen as functions $\Omega\rightarrow\mathbb{R}$. We then have:
\begin{equation}\label{eq:partialvx}
    \partial_xv(\xi)=\partial_\theta v(\xi)\partial_x\theta(\xi)+\partial_\sigma v(\xi)\partial_x\sigma(\xi)+\partial_\rho v(\xi)\partial_x\rho(\xi),\quad\xi\in\Omega.
\end{equation}
and
\begin{equation}\label{eq:partialvy}
    \partial_y v(\xi)=\partial_\theta v(\xi)\partial_y\theta(\xi)+\partial_\sigma v(\xi)\partial_y\sigma(\xi)+\partial_\rho v(\xi)\partial_y\rho(\xi),\quad\xi\in\Omega.
\end{equation}

If $\xi\in L_1$, Condition \textbf{(A)} then requires both $\partial_x v(\xi)$ and $\partial_y v(\xi)$ to be zero. The restriction of $v$ to $\Sigma$ is given by $f$, and so the tangential derivatives of both functions coincide at $\xi$, i.e., $\partial_\theta v(\xi)=\partial_\theta f(\xi)$ and $\partial_\sigma v(\xi)=\partial_\sigma f(\xi)$. The normal derivative of $v$ at $\xi$ is prescribed by $g$, that is, we have the equality $\partial_\rho v(\xi)=g(\xi)$. On $L_1$, the gradients of the coordinates $\theta$, $\sigma$ and $\rho$ can be identified component-wise with the vectors $\mathbf{t}$, $\mathbf{k}$ and $\mathbf{n}$ of the moving frame $\mathcal{B}$. 

Denoting the components of each of this vectors with a subindex, the constraints on the two data $f$ and $g$ resulting from Condition \textbf{(A)} are then:
\begin{equation}\label{eq:constraintA}
\textbf{(A):}
    \begin{cases}
        \partial_x v(\xi)=\partial_\theta f(\xi)t_1(\xi)+ \partial_\sigma f(\xi)k_1(\xi)+ g(\xi)n_1(\xi)=0,\\
        \partial_y v(\xi)=\partial_\theta f(\xi)t_2(\xi)+ \partial_\sigma f(\xi)k_2(\xi)+ g(\xi)n_2(\xi)=0,\quad\forall\xi\in L_1.
    \end{cases}
\end{equation}
So, Condition \textbf{(A)} leads to a homogeneous linear system of two equations and three unknowns, that links the derivatives of $f$ and the normal datum $g$.

Setting the auxiliary function:
\begin{align}\label{eq:auxmu}     \mu(\xi):=n_1(\xi)^2+n_2(\xi)^2,
\end{align}
since $\mathbf{k}$ satisfies \eqref{eq:condk}, the two components $n_1$ and $n_2$ cannot simultaneously vanish, and consequently $\mu(\xi)>0$ for all $\xi\in L_1$. We then easily get
\begin{equation}
    g(\xi)=-\frac{1}{\mu(\xi)}\Big(\partial_\theta f\big(\xi)(t_1n_1+t_2n_2\big)+\partial_\sigma f(\xi)\big(k_1n_1+k_2n_2\big)\Big),\quad\forall\xi\in L_1.
\end{equation}
By recalling that the vectors in $\mathcal{B}$ form an orthonormal basis, we can finally write the constraint that \textbf{(A)} imposes on $g$ as:
\begin{equation}\label{eq:constraintA1}
    \textbf{(A1):}\quad g(\xi)=\frac{1}{\mu(\xi)}n_3(\xi)\big(\partial_\theta f(\xi)t_3(\xi)+\partial_\sigma f(\xi)k_3(\xi)\big),\quad\forall\xi\in L_1.
\end{equation}
That is, Condition \textbf{(A)} completely fixes the values of $g$ along $L_1$ in terms of those of the tangential derivatives of $f$ and the components of the vectors in the moving frame $\mathcal{B}$.

Next, let us see how  Condition \textbf{(A)} links the two tangential derivatives $\partial_\theta f$ and $\partial_\sigma f$. For this, we substitute Equation \eqref{eq:constraintA1} into Equation \eqref{eq:constraintA}. After straightforward manipulations using the orthogonality relations between the vectors of $\mathcal{B}$, the first equality of the system becomes:
\begin{equation}\label{eq:step1constraintA2}
    n_2(\xi)\big(\partial_\sigma f(\xi)t_3(\xi)-\partial_\theta f(\xi)k_3(\xi)\big)=0,\quad\xi\in L_1.    
\end{equation}
Likewise, for the second equality coming from \textbf{(A)} in the system \eqref{eq:constraintA}, we obtain:
\begin{equation}\label{eq:step2constraintA2}
    n_1(\xi)\big(\partial_\theta f(\xi)k_3(\xi)-\partial_\sigma f(\xi)t_3(\xi)\big)=0,\quad\xi\in L_1.    
\end{equation}
The terms in between brackets in \eqref{eq:step1constraintA2} and \eqref{eq:step2constraintA2} only differ in sign. Since $n_1(\xi)$ and $n_2(\xi)$ cannot be zero at the same time, it follows that the constraint on $f$ resulting from \textbf{(A)} is:
\begin{equation}\label{eq:constraintA2}
    \textbf{(A2):}\quad \partial_\sigma f(\xi)t_3(\xi)-\partial_\theta f(\xi)k_3(\xi)=0,\quad\forall\xi\in L_1.
\end{equation}
It is clear that Condition \textbf{(A)} is equivalent to Conditions \textbf{(A1)} and \textbf{(A2)}.

\subsubsection{Step 2: Choice of a general ansatz for $f$ and $g$.}\label{subs:existence}

Since $t_3(\xi)$ is zero at the horizontal points of $L_1$, it is not possible to divide Equation \eqref{eq:constraintA2} by it. Thus, for \textbf{(A2)} to be satisfied, the two tangential derivatives of $f$ have to be of the form
\begin{equation}\label{eq:fderivatives}
    \begin{cases}
        \partial_\theta f(\xi)=\beta(\xi)t_3(\xi),\\
        \partial_\sigma f(\xi)=\beta(\xi)k_3(\xi),\quad{\xi\in L_1},
    \end{cases}
\end{equation}
for some arbitrary analytic function $\beta:\Omega\rightarrow\mathbb{R}$.

The result of substituting \eqref{eq:fderivatives} into \eqref{eq:constraintA1} is:
\begin{equation}
    g(\xi)=\frac{1}{\mu(\xi)}\beta(\xi)n_3(\xi)\big(t_3(\xi)^2+k_3(\xi)^2\big),\quad\xi\in L_1.
\end{equation}
Because of the orthonormality of the vectors in $\mathcal{B}$, the term appearing in between brackets here is equal to $n_1(\xi)^2+n_2(\xi)^2$, that is, the function $\mu(\xi)$. It follows that the constraint is:
\begin{equation}\label{eq:gsimplified}
    g(\xi)=\beta(\xi)n_3(\xi),\quad\forall\xi\in L_1.
\end{equation}
So, Condition \textbf{(A)} requires that, for all $\xi\in L_1$, the triple of first-order data $(\partial_\theta f(\xi),\partial_\sigma f(\xi),g(\xi))$ has to be proportional to the projection of $\mathbf{e_3}$ onto the basis of the tangent space $T_{\xi}\Omega$ given by \eqref{eq:omegatangentbasis}.
\begin{remark}
The function $\beta$ is related to the derivative $\partial_t v$. Indeed, using the chain rule we get
\begin{equation}
       \partial_t v(\xi)=\partial_\theta f(\xi)t_3(\xi)+\partial_\sigma f(\xi)k_3(\xi)+g(\xi)n_3(\xi),\quad\xi\in L_1.
\end{equation}
Now, we substitute into this expression Equations \eqref{eq:fderivatives} and \eqref{eq:gsimplified}, yielding
\begin{equation}\label{eq:finalvt}
        \partial_t v(\xi)=\beta(\xi)\big(t_3(\xi)^2+k_3(\xi)^2+n_3(\xi)^2\big)
        =\beta(\xi),\quad\forall\xi\in L_1.
\end{equation}
Here we have used that $t_3^2+k_3^2+n_3^2\equiv 1$ thanks to the orthonormality of the vectors involved in $\mathcal{B}$. Thus, the time derivative $\partial_t v$ coincides with $\beta$ in $L_1$, and Constraint \textbf{(A)} leads to $\nabla_3 v(\xi)\equiv\beta(\xi)\mathbf{e_3}$ for $\xi\in L_1$.
\end{remark}

In order to show that there exists an analytic function $f$ whose derivatives are given by \eqref{eq:fderivatives}, we make use of the following general lemma. Its proof is relegated to Appendix \ref{app:appendixA}.

\begin{lemma}\label{lemma:existence}
    Let $\Lambda$ be an analytic curve in $\mathbb{R}^3$, and suppose that $\mathbf{r_\Lambda}:\Lambda\rightarrow\mathbb{R}^3$ is a $C^\omega$ vector field transverse to it. Then, there exists an analytic function $h_\Lambda:\Omega\rightarrow\mathbb{R}$, defined in an open neighborhood $\Omega$ of $\Lambda$ in $\mathbb{R}^3$, whose gradient's restriction $\nabla_3 h_\Lambda|\Lambda$ is equal to $\mathbf{r_\Lambda}$.
\end{lemma}

So, for our case, consider the vector field:
\begin{align}
        \zeta(\xi):=\beta(\xi)t_3(\xi)\mathbf{t}(\xi)+\beta(\xi)k_3(\xi)\mathbf{k}(\xi)+\mathbf{n}(\xi),\quad\xi\in L_1.
\end{align}
It is analytic, and transverse to $L_1$ thanks to the third term $\mathbf{n}(\mathbf{x})$. Lemma \ref{lemma:existence} then implies that there exists a function $f_1:\Omega\rightarrow\mathbb{R}$ defined in an open neighborhood $\Omega\subset\mathbb{R}^3$ of $L_1$ such that $\nabla_3 f_1|_{L_1}=\zeta$. If necessary we reduce $\sigma_0$, $\sigma_1$ and $\rho_1$ in \eqref{eq:sigmaparam} and \eqref{eq:omegaparam} so that $\Omega$ is the same as before. Of course, $f_1$ is an analytic function whose derivatives are given by \eqref{eq:fderivatives}, as desired.

However, it is convenient to modify the function $f_1$ so that we can prescribe the value of  $\partial_{\sigma\sigma} f$ on $L_1$ all while still satisfying Equation~\eqref{eq:fderivatives}. If $F:\Omega\rightarrow\mathbb{R}$ is an arbitrary $C^\omega$ function, we define the analytic function:
\begin{align}
        f(\xi):= f_1(\xi)-\frac{1}{2}\sigma(\xi)^2\Big(\partial_{\sigma\sigma}f_1(\xi)-F(\xi)\Big).
\end{align}
It is clear that $f=f_1$ on $L_1$ and, moreover, $f$ still satisfies Equation~\eqref{eq:fderivatives}. Additionally, a straightforward computation shows that
\begin{equation}
    \partial_{\sigma\sigma}  f(\xi)=F(\xi)
\end{equation}
if $\xi\in L_1$, so that we can prescribe this second-order partial derivative on $L_1$.

\begin{remark}
Using Equation~\eqref{eq:fderivatives} we immediately get
\begin{equation}\label{eq:partialfthetatheta}
    \partial_{\theta\theta} f(\xi)=\partial_\theta\beta(\xi)t_3(\xi)+\beta(\xi)\partial_\theta t_3(\xi),\quad\xi\in L_1.
\end{equation}
and
\begin{equation}\label{eq:partialfsigmatheta}
    \partial_{\sigma\theta} f(\xi)=\partial_\theta\beta(\xi)k_3(\xi)+\beta(\xi)\partial_\theta k_3(\xi),\quad\xi\in L_1.
\end{equation}
\end{remark}

Concerning the function $g$, the only constraint is that its values along $L_1$ have to be given by \eqref{eq:gsimplified}. So, assume $G:\Omega\rightarrow\mathbb{R}$ is an arbitrary $C^\omega$ function, and define:
\begin{align}\label{eq:gfinal}
        g(\xi)=\beta(\xi) (\mathbf{v}_{\rho,\Omega})_3(\xi)+\sigma(\xi) G(\xi).
\end{align}
Here, $(\mathbf{v}_{\rho,\Omega})_3$ is the third component of the vector field $\mathbf{v}_{\rho,\Omega}$. Clearly, this is an analytic function in $\Omega$ and it satisfies \eqref{eq:gsimplified}. 

In what follows, we shall take the two analytic functions $f$ and $g$ defined before, which automatically satisfy Condition~\textbf{(A)}, and we shall exploit the freedom of the functions $\beta,F,G$ in order to satisfy conditions \textbf{(B)} and \textbf{(C)}.

\subsubsection{Step~3: Final adjustments to the Cauchy data }\label{subs:stability}

In order to simplify the notation, from now on we use $v_x$ to denote $\partial_x v$ and $v_y$ to do so for $\partial_y v$. We first show that the three arbitrary (analytic) functions $\beta,F,G$ can be chosen so that Condition \textbf{(C)} is verified:
\begin{equation}
    \rank\big(\nabla_3 v_x(\xi),\nabla_3 v_y(\xi)\big)=2,\quad\forall\xi\in L_1.
\end{equation}
To this end, it is convenient to use the coordinate system $(\theta,\sigma, \rho)$. Since $v_x(\xi)=v_y(\xi)=0$ for all $\xi\in L_1$, then $\partial_\theta v_x=\partial_\theta v_y=0$ on $L_1$. Accordingly, the only nontrivial part of the Hessian matrix of $v$ at $\xi\in L_1$ is its projection onto the normal space $N_{\xi}L_1$. This way, $(\nabla_3 v_x(\xi),\nabla_3 v_y(\xi))$ is equivalent, under an orthogonal change of basis, to a $2\times 2$ matrix whose components are the projections of $\nabla_3 v_x$ and $\nabla_3 v_y$ along $\mathbf{k}(\xi)$ and $\mathbf{n}(\xi)$, i.e., the derivatives of $v_x$ and $v_y$ with respect to the coordinates $\sigma$ and $\rho$ on $L_1$. We thus define the \textit{normal stability matrix}:
\begin{equation}\label{eq:normalmatrix}
    H_N(\xi):=
    \begin{pmatrix}
        \partial_\sigma v_x(\xi) & \partial_\rho v_x(\xi)\\
        \partial_\sigma v_y(\xi) & \partial_\rho v_y(\xi)
    \end{pmatrix},\quad\xi\in L_1.
\end{equation}
We claim that $v$ satisfies Condition \textbf{(C)} if and only if the determinant of $H_N(\xi)$ is nonzero for all $\xi\in L_1$:
\begin{align}\label{eq:normalstabilitydet}
\cD_N(\xi):=\partial_\sigma v_x(\xi)\partial_\rho v_y(\xi)-\partial_\rho v_x(\xi)\partial_\sigma v_y(\xi)\neq 0
\end{align}
for all $\xi\in L_1$. We say that $\cD_N$ is the \textit{normal stability determinant}. 

Indeed, observe that the type of the nondegenerate (spatial) critical point $\xi\in L_1$ of $v$, can be determined from $\Delta v(\xi)=\beta(\xi)$, $\xi\in L_1$, cf. Equation~\eqref{eq:finalvt}, and the determinant:
\begin{align}
        \cD_{xy}(\xi):=v_{xx}(\xi)v_{yy}(\xi)-v_{xy}(\xi)v_{yx}(\xi).
\end{align}
The point $\xi\in L_1$ is an extremum if $\cD_{xy}(\xi)>0$ and a saddle if $\cD_{xy}(\xi)<0$, where in the first case it is a maximum if $\beta(\xi)<0$ and a minimum if $\beta(\xi)>0$.

An elementary computation using the chain rule allows us to compute $\cD_{xy}(\xi)$ in terms of the normal stability determinant:
\begin{align}\label{eq:deltaXY}
        \cD_{xy}(\xi)=\cD_N(\xi)t_3(\xi),\quad{\xi\in L_1}.
\end{align}
Analogously, setting
$$
\cD_{xt}(\xi):= v_{xx}(\xi)v_{tt}(\xi)-v_{xt}(\xi)v_{tx}(\xi),
$$
and
$$
\cD_{yt}(\xi):= v_{yy}(\xi)v_{tt}(\xi)-v_{yt}(\xi)v_{ty}(\xi),
$$
we easily find
\begin{align}
        \cD_{xt}(\xi)=-\cD_N(\xi)t_2(\xi),\quad{\xi\in L_1},
\end{align}
and
\begin{align}
        \cD_{yt}(\xi)=\cD_N(\xi)t_1(\xi),\quad{\xi\in L_1}.
\end{align}
We conclude that the three minors $\cD_{xy}$, $\cD_{xt}$ and $\cD_{yt}$ vanish simultaneously at $\xi\in L_1$ if and only the normal stability determinant $\cD_N(\xi)$ does, which shows that Condition \textbf{(C)} is fulfilled if and only if $\cD_N$ does not have zeroes in $L_1$, as we wanted to show. 

Our next goal is to find a formula for $\cD_N$ in terms of the functions $\beta,F,G$. We first give an expression for $\partial_\sigma v_x$ and $\partial_\sigma v_y$. An easy computation using the data $f$ and $g$, as in Equation~\eqref{eq:constraintA}, and the vectors of $\mathcal{B}$ for the first-order derivatives of the coordinates, yields:
\begin{align}\label{eq:vxomega}
    \begin{split}
        \partial_\sigma v_x(\xi)&=f_{\theta\sigma}(\xi)t_1(\xi)+
        f_{\sigma\sigma}(\xi)k_1(\xi)+g_\sigma(\xi)n_1(\xi)\\
        &+f_\theta(\xi)\partial_\sigma\theta_x(\xi)+f_\sigma(\xi)\partial_\sigma\sigma_x(\xi) + g(\xi)\partial_\sigma\rho_x(\xi),\quad\xi\in L_1.
    \end{split}
\end{align}
Similarly, for $v_y$ we obtain:
\begin{align}\label{eq:vyomega}
    \begin{split}
        \partial_\sigma v_y(\xi)&=f_{\theta\sigma}(\xi)t_2(\xi)+
        f_{\sigma\sigma}(\xi)k_2(\xi)+g_\sigma(\xi)n_2(\xi)\\
        &+f_\theta(\xi)\partial_\sigma\theta_y(\xi)+f_\sigma(\xi)\partial_\sigma\sigma_y(\xi) + g(\xi)\partial_\sigma\rho_y(\xi),\quad\xi\in L_1.
    \end{split}
\end{align}
The derivatives of the Cauchy data $f$ and $g$ in these last two expressions have been described in Subsection \ref{subs:existence}. Specifically, $\partial_\theta f$ and $\partial_\sigma f$ appear in \eqref{eq:fderivatives}, $\partial_{\theta\sigma} f$ in \eqref{eq:partialfsigmatheta}, $\partial_{\sigma\sigma} f$ is the arbitrary function $F$, and $\partial_\sigma g$ is $G$.

Regarding the two partial derivatives $\partial_\rho v_x$ and $\partial_\rho v_y$, it is convenient to write a system in terms of the two second-order derivatives which we do have direct access to, $\partial_\sigma v_x$ and $\partial_\sigma v_y$, making all the data appear directly. To this end, we write  $v_{xx}$ and $v_{yy}$ using the chain rule:
\begin{align}\label{eq:interstepB1}
    v_{xx}(\xi)&=\partial_\theta v_x(\xi) \partial_x \theta(\xi) + \partial_\sigma v_x(\xi) \partial_x \sigma(\xi) +\partial_\rho v_x(\xi) \partial_x \rho(\xi),\quad\xi\in\Omega,\\
    v_{yy}(\xi)&=\partial_\theta v_y(\xi) \partial_y \theta(\xi) + \partial_\sigma v_y(\xi) \partial_y \sigma(\xi) +\partial_\rho v_y(\xi) \partial_y \rho(\xi),\quad\xi\in\Omega.
\end{align}

If $\xi$ belongs to $L_1$, then both $\partial_\theta v_x(\xi)$ and $\partial_\theta v_y(\xi)$ are zero. Moreover, from \eqref{eq:finalvt}, it also is $v_t(\xi)=\beta(\xi)$. Introducing all this back into the heat equation, along with the components of the basis vectors $\mathbf{k}(\xi)$ and $\mathbf{n}(\xi)$ for the partial derivatives of the coordinate functions $\sigma$ and $\rho$, we get the first constraint:
\begin{equation}\label{eq:normaleq1}
    \quad \partial_\sigma v_x(\xi) k_1(\xi) + \partial_\rho v_x(\xi) n_1(\xi) + \partial_\sigma v_y(\xi) k_2(\xi) + \partial_\rho v_y(\xi) n_2(\xi) =\beta(\xi),\quad\forall \xi\in L_1.
\end{equation}
The second equation from where to solve for $\partial_\sigma v_x$ and $\partial_\sigma v_y$ is simply $v_{xy}=v_{yx}$. Proceeding as before, we find the second constraint:
\begin{equation}\label{eq:normaleq2}
    \partial_\sigma v_x(\xi)k_2(\xi)  + \partial_\rho v_x(\xi)n_2(\xi)  = \partial_\sigma v_y(\xi)k_1(\xi)  + \partial_\rho v_y(\xi)n_1(\xi) ,\quad\forall \xi\in L_1.
\end{equation}

Equations \eqref{eq:normaleq1} and \eqref{eq:normaleq2} form a linear system for the two unknowns $\partial_\rho v_x$ and $\partial_\rho v_y$. Writing both together, in matrix notation, and omitting the argument $\xi$, we obtain the unique solution
\begin{equation}\label{eq:finalnormal}
    \begin{pmatrix}
        \partial_\rho v_x\\
        \partial_\rho v_y
    \end{pmatrix}
     =
     \frac{1}{\mu}\left[
     \begin{pmatrix}
        n_1 & n_2\\
        n_2 & -n_1
    \end{pmatrix}
     \begin{pmatrix}
         -k_1 & -k_2\\
         -k_2 & k_1
     \end{pmatrix}
     \begin{pmatrix}
         \partial_\sigma v_x\\
         \partial_\sigma v_x
     \end{pmatrix}
     +
     \beta
     \begin{pmatrix}
         n_1\\
         n_2
     \end{pmatrix}
     \right]
     ,\quad\xi\in L_1.
\end{equation}
Recall that the the function $\mu$ does not vanish on $L_1$. This provides the expression we wanted for $\partial_\rho v_x(\xi)$ and $\partial_\rho v_y(\xi)$ in $L_1$. Indeed, all the terms on the right side are known; they depend only on the Cauchy data $f$ and $g$, their tangential derivatives, and the noncharacteristic surface $\Sigma$.

Then, if we introduce \eqref{eq:vxomega} and \eqref{eq:vyomega}, together with \eqref{eq:finalnormal}, into the formula for $\cD_N$ appearing in \eqref{eq:normalstabilitydet}, we obtain an  expression for the normal stability determinant:
\begin{equation}\label{eq:finalnormaldeterminant}
    \cD_N(\xi)=\frac{1}{\mu(\xi)}\Big( \beta(\xi)\big( \partial_\sigma v_x(\xi) n_2(\xi)-\partial_\sigma v_y(\xi) n_1(\xi)\big)-t_3(\xi)\big(\partial_\sigma v_x(\xi)^2+\partial_\sigma v_y(\xi)^2\big)\Big),\;\xi\in L_1.
\end{equation}
It is obvious from this expression that for Condition \textbf{(C)} to be satisfied, then $\beta$ is required to be nonzero at the height extrema of $L_1$. More generally, either $\cD_N(\xi)>0$ or $\cD_N(\xi)<0$ for all $\xi\in L_1$.

In parallel, Condition \textbf{(B)} constraints which type of critical point is each $\xi\in L_1$. For our choice of positive and negative intervals in Theorem \ref{thm:maintheoremR2}, the critical points of $v$ for which $t_3(\xi)$ is positive have to be extrema, and those for which $t_3(\xi)$ is negative, saddles. So, looking at \eqref{eq:normalstabilitydet} and \eqref{eq:deltaXY}, in order to satisfy Condition~\textbf{(B)}, it has to be 
$$\cD_N(\xi)>0$$  
for all $\xi\in L_1$. Setting the \textit{reduced normal determinant} $\eta(\xi):=\mu(\xi)\cD_N(\xi)$, we infer from the previous discussion that both Conditions~\textbf{(B)} and~\textbf{(C)} hold if and only if 
$$\eta(\xi)>0,\quad\forall \xi\in L_1.$$

We claim that we can choose the free analytic functions $F$, $G$ and $\beta$ so this is indeed verified, as stated next.

\begin{lemma}\label{lemma:FGbeta}
    There exist three analytic functions $\beta, F,G:\Omega\rightarrow\mathbb{R}$ such that $\eta(\xi)>0$ for all $\xi\in L_1$. Furthermore, $\beta$ can be chosen so that the critical points of $v$ in $L_1^+$ (extrema) are in accordance with our initial choice of maxima and minima.
\end{lemma}

The proof of Lemma \ref{lemma:FGbeta} is presented in Appendix \ref{app:appendixB}. In view of this result, if we take the functions $F,G,\beta$ provided by Lemma~\ref{lemma:FGbeta}, from the construction introduced in the previous subsection we obtain a pair of Cauchy data $f$ and $g$ with all the desired properties,  which we can directly introduce in the Cauchy problem \eqref{eq:heatcauchy}. Then, as stated at the beginning, the existence of the heat solution $v$ is a consequence of the Cauchy-Kovalevskaya theorem, and by construction of the Cauchy data it satisfies the Conditions~\textbf{(A)}, \textbf{(B)} and~\textbf{(C)}. By reducing $\Omega$ if necessary, the statement of Proposition \ref{prop:heat} then follows.

\section{Proof of Proposition~\ref{prop:approx}}\label{section:approx}

First, we observe that the spacetime domain $\Omega$ satisfies the topological condition that the complement of its intersection with each $t$-hyperplane for $t>0$ is connected. Then, a direct application of  \cite[Theorem 1.2.]{enciso_approximation_2019} yields the existence of a function $u_0\in C^\infty_c(\mathbb{R}^2)$ such that $u:=e^{t\Delta}u_0$ is a global solution to the heat equation satisfying:
\begin{equation}\label{eq:ineq1}
    \|v-u\|_{C^{r+1}(\Omega)}<\delta_1.
\end{equation} 
As usual, $u:=e^{t\Delta}u_0:\mathbb{R}^2\times\mathbb{R}_+\rightarrow\mathbb{R}$ is the unique solution to the Cauchy problem
\begin{equation}
    \begin{cases}
       \partial_t u= \Delta u,\\
        u(x,y,0)=u_0(x,y),\quad\text{with}\quad(x,y,t)\in\mathbb{R}^2\times\mathbb{R}_+,
    \end{cases}
\end{equation}
such that $u(x,y,t)\rightarrow0$ when $|(x,y)|\rightarrow\infty$ for all $t>0$. 

Next, let $\delta_2>0$ be any real number, and consider the following Cauchy problem:
\begin{equation}\label{eq:cauchystreamscaled}
    \begin{cases}
        \partial_t \tilde{\psi}+\delta_2\Delta^{-1}\big(\nabla\tilde{\psi}\cdot\nabla^\perp(\Delta\tilde{\psi})\big)=\Delta\tilde{\psi},\\
        \tilde{\psi}(x,y,0)=u_0(x,y),\quad\text{with}\quad(x,y,t)\in\mathbb{R}^2\times\mathbb{R}_+.
    \end{cases}
\end{equation}
By {Duhamel's Principle} \cite[Proposition 1.35.]{tao_nonlinear_2006}, it is equivalent to the integral equation:
\begin{equation}\label{eq:duhamel}
    \tilde{\psi}(x,y,t)=e^{t\Delta}u_0 - \delta_2\int_0^t e^{(t-s)\Delta}\Delta^{-1}(\nabla\tilde{\psi}\cdot\nabla^\perp(\Delta\tilde{\psi}))(x,y,s)ds.
\end{equation}
The first term on the right side of the equality is the global heat solution $u$ and the second one is bounded for finite times $t$. If the parameter $\delta_2$ is small enough, then there exists a global solution $\tilde{\psi}$ such that 
\begin{equation}\label{eq:ineq2}
    \|u-\tilde{\psi}\|_{C^{r+1}(\mathbb{R}^2\times[-T,T])}<\delta_2C_T,
\end{equation}
for any time $T>0$. Here $C_T$ is a constant which only depends on $T$, but not on the initial condition $u_0$ or on $\delta_2$ itself.

Finally, take a big enough $T$ so that $\overline{\Omega}\subset\mathbb{R}^2\times[0,T]$. Then, by \eqref{eq:ineq2}:
\begin{equation}
    \|u-\tilde{\psi}\|_{C^{r+1}(\Omega)}<\delta_2C_T,
\end{equation}
and using~\eqref{eq:ineq1}, we obtain:
\begin{equation}\label{eq:finalbound}
    \|v-\tilde{\psi}\|_{C^{r+1}(\Omega)}\leq  \delta_1+C_T\delta_2.
\end{equation}
Taking small enough constants $\delta_1$ and $\delta_2$ the proposition follows.

\section{The Case of the Navier-Stokes Equations on the Torus}\label{section:torus}

\noindent In this section we prove that there are analogous statements to Theorem \ref{thm:mergingsplittingR2} and \ref{thm:maintheoremR2} for the case of the torus. The only caveat is that we have to allow the link $L$ to undergo a {scaling}. The reason is that the global approximation theorems in \cite{enciso_approximation_2019} for local heat solutions do not apply in $\mathbb{T}^2$ because of its compactness. To overcome this difficulty, the idea is to look at small scales, for times of order $\varepsilon$ and spatial variables of order $\sqrt{\varepsilon}$, where we can apply~\cite[Lemma 9.1.]{enciso_approximation_2019}.

Let $B_1$ be the image in $\mathbb{T}^2$ of a unit ball in the square $[0,2\pi]\times[0,2\pi]$. Choose a real number $0<\eta<1$, and consider the linear map:
\begin{align}\label{eq:scalingeta}
    \begin{split}
        \Lambda_\eta:&B_1\times\mathbb{R}_+\rightarrow B_1\times\mathbb{R}_+\\
        &(x,y,t)\mapsto(\sqrt{\eta}x,\sqrt{\eta}y,\eta t).
    \end{split}
\end{align}

Let $L$ be a link, possibly open, in $B_1\times\mathbb{R}_+\subset \mathbb{T}^2\times\mathbb{R}_+$. Assume that $L$ is partitioned into $L^+$ and $L^-$ as in Section \ref{section:intro} through the choice of an orientation $\zeta$. Pick a subfamily $L^+_{\max}$ of positive intervals in $L^+$, and define $L^+_{\min}:=L^+\backslash L^+_{\max}$.

Our main result for the case of the Navier-Stokes equations on $\mathbb{T}^2$ (analogous to Theorem \ref{thm:maintheoremR2}) is:

\begin{theorem}\label{thm:maintheoremT2}
    Fix a positive integer $r\geq 1$ and let $\varepsilon>0$ be any real number. Let $L$ be an arbitrary link in $B_1\times\mathbb{R}_+$, possibly open, with an orientation $\zeta$. There exists a solution $\psi:\mathbb{T}^2\times\mathbb{R}_+\rightarrow\mathbb{R}$ to the Navier-Stokes equations in $\mathbb{T}^2$, a tubular neighborhood $\Omega$ of $L$, a real number $0<\eta\ll1$, and a diffeomorphism $\Phi:\Lambda_\eta(\Omega)\rightarrow\Lambda_\eta(\Omega)$ such that:
    \begin{itemize}
        \item $\|\Phi-\mathbf{id}\|_{C^r(\Omega^\prime)}<\varepsilon$, with $\Omega^\prime:=\Lambda_\eta(\Omega)$.
        \item The set of critical points of $\psi$ in $\Omega^\prime$ is exactly $L_2:=(\Phi\circ\Lambda_\eta)(L)$.
    \end{itemize}
    Moreover, the type of each critical point of $\psi$ along $L_2$ is known:
    \begin{itemize}
        \item Consider the orientation induced by $\zeta$ in $L_2$ via the map $\Phi\circ\Lambda_\eta$. Then, the critical points of $\psi$ in the positive intervals $L_2^+$ are extrema, and those in $L_2^-$ are saddles.
        \item Consider the decomposition of $L_2^+$ into $L_{2,\max}^+$ and $L_{2,\min}^+$ coming from $L^+$. The critical points of $\psi$ in the former intervals are maxima, and those in the latter are minima.
    \end{itemize}
\end{theorem}

\noindent As in $\mathbb{R}^2$, the situation here is structurally stable. That is, any stream function that is close enough to $\psi$ in the $C^r$ norm realizes the scaled link $\Lambda_\eta(L)$ as an isolated subset of its critical set up to a diffeomorphism. Remark \ref{remarks:commentstheoremR2} is valid here too.

At each horizontal point of the new link $L_2$ there is a splitting or merging of critical points for the stream function $\psi$. It follows that Theorem \ref{thm:mergingsplittingR2} is still valid on $\mathbb{T}^2$, just having to allow the merging and splitting arcs $\Gamma$ to undergo first the scaling given by $\Lambda_\eta$. {The arcs are defined here as in Section \ref{section:intro}.} The following is a direct corollary of Theorem~\ref{thm:maintheoremT2}.

\begin{theorem}\label{thm:mergingsplittingT2}
    Let us fix a positive integer $r\geq 1$ and some small $\varepsilon>0$. Consider $k$ connecting arcs {$\Gamma_i$ in $B_1\times \mathbb{R}_+$}, which we can assume to be merging for $1\leq i\leq k'$ and splitting for $k'<i\leq k$. 
    
    Then there is some initial datum {$\psi_0\in C^\infty_c(\mathbb{T}^2)$} for which the associated solution {$\psi:\mathbb{T}^2\times\mathbb{R}_+\to\mathbb{R} $} to the Navier--Stokes equations realizes that bifurcation pattern, up to {a scaling and} a small deformation. More precisely, there exist {a real number $0<\eta\ll1 $ and a diffeomorphism $\Phi$ of~$\mathbb{T}^2\times\mathbb{R}_+$ with 
$\|\Phi-\mathbf{id}\|_{C^r(\mathbb{T}^2\times\mathbb{R}_+)}<\varepsilon$} such that the spacetime curves {$\widetilde{\Gamma}_i:=(\Phi\circ\Lambda_\eta)(\Gamma_i)$} consist of isolated stagnation points of the fluid. %\textcolor{blue}{That is, $\nabla \psi\circ \widetilde{\Gamma}_i\equiv 0$ and there exists an open tubular neighborhood $\Omega_i\subset B_1\times\mathbb{R}_+$ whose subset of critical points of $\psi$ is the trace of $\tilde{\Gamma}_i$.}

Furthermore, the functions $h\circ\widetilde{\Gamma_i}:[0,1]\to\mathbb{R}$ are still strictly concave or convex, and attain their maximum or minimum at a unique interior point $s_i\in(0,1)$. The splitting (for $i>k'$) or merging (for $i\leq k'$) of stagnation points happens at the spacetime point $\widetilde{\Gamma}_i(s_i)$, which is a degenerate critical point of~$\psi$. The stagnation point $\widetilde{\Gamma}_i(s)$ can be chosen to be a nondegenerate local extremum for $s\in[0,s_i)$ and a nondegenerate saddle for $s\in(s_i,1]$.
\end{theorem}

Let us sketch the proof of Theorem~\ref{thm:maintheoremT2}. The first observation is that Proposition \ref{prop:heat} is still valid in $B_1\times\mathbb{R}_+$. Accordingly, let $v:\Omega\rightarrow\mathbb{R}$ be a local heat solution with properties \textbf{(A)}, \textbf{(B)} and \textbf{(C)}. The domain $\Omega\subset B_1\times\mathbb{R}_+$ is a neighborhood of $L_1:=\Phi_1(L)$, where $\Phi_1:\Omega\rightarrow\Omega$ satisfies $\|\Phi_1-\mathbf{id}\|_{C^r(\Omega)}<\varepsilon_1$; by reducing $\Omega$ if needed, we can safely assume that its intersection with each time-$t$ slice has a connected complement in $B_1$. Then, for any $\varepsilon_2>0$, by \cite[Lemma 9.1.]{enciso_approximation_2019}, there exists a real number {$0<\eta\ll1$} and a global heat solution on the torus $u:\mathbb{T}^2\times\mathbb{R}_+\rightarrow\mathbb{R}$, such that:
\begin{equation}
    \|u\circ\Lambda_\eta-v\|_{C^{r+1}(\Omega)}<\varepsilon_2.
\end{equation}
Passing to Navier-Stokes stream function as in Section \ref{section:proofmain}, there exists a diffeomorphism $\Phi_2:\Omega^\prime\rightarrow\Omega^\prime$ such that $\|\Phi_2-\mathbf{id}\|_{C^r(\Omega^\prime)}<\varepsilon_2$ and $\Phi_2(\Lambda_\eta(L_1))$ is an isolated subset of critical points of a solution $\psi:\mathbb{T}^2\times\mathbb{R}_+\rightarrow\mathbb{R}$ to the Navier-Stokes equations, where $\Omega^\prime:=\Lambda_\eta(\Omega)$. Now, we just define $\Phi:=\Phi_2\circ(\Lambda_\eta\circ \Phi_1\circ \Lambda_\eta^{-1}):\Lambda_\eta(\Omega)\rightarrow\Lambda_\eta(\Omega)$, and since $\varepsilon_1$ and $\varepsilon_2$ are arbitrary, the result immediately follows.

\section*{Acknowledgments} This work has received funding from the European Research Council (ERC) under the European Union's Horizon 2020 research and innovation programme through the grant agreement~862342 (A.E.). I.B.\ is supported by MCIN under the scholarship PREP2022-000014 at CSIC. The authors are partially supported by the grants CEX2023-001347-S, RED2022-134301-T and PID2022-136795NB-I00 (I.B., A.E., and D.P.-S.) funded by MCIN/AEI/10.13039/501100011033, and Ayudas Fundaci\'on BBVA a Proyectos de Investigaci\'on Cient\'ifica 2021 (D.P.-S.).  

\appendix
\section{Proof of Lemma~\ref{lemma:existence}}\label{app:appendixA}

In view of {Cartan's Analytic Extension Theorem}~\cite[Theorem 7.4.8]{hormander_introduction_1990}, we can safely assume that $\mathbf{r_\Lambda}$ is an analytic vector field defined on the whole $\mathbb R^3$.
Next, let $\mathbf{t_\Lambda}:\Lambda\rightarrow\mathbb{R}^3$ be a unit tangent vector field along the curve $\Lambda$ (for some prescribed orientation), and define the vector field
$$
\mathbf{k_{\Lambda}}(\xi):=\frac{\mathbf{t_\Lambda}(\xi)\times\mathbf{r_\Lambda}(\xi)}{\big\|\mathbf{t_\Lambda}(\xi)\times\mathbf{r_\Lambda}(\xi)\big\|}
$$
for $\xi\in \Lambda$. By construction, the vector field $\mathbf{k_\Lambda}$ is analytic and transverse to $\Lambda$.

We denote by $\Sigma_\Lambda$ a ruled surface obtained by sliding the points of $\Lambda$ along the direction prescribed by $\mathbf{k_{\Lambda}}$; this surface is clearly analytic and contains the curve $\Lambda$. Now, let $\mathbf{N_{\Lambda}}:\Sigma_\Lambda\rightarrow\mathbb{R}^3$ be a normal unit vector field on $\Sigma_\Lambda$. It is clear from the construction that $\mathbf{N_{\Lambda}}|_{\Lambda}$ is proportional to $\mathbf{r_\Lambda}$, so we can assume without any loss of generality that $\mathbf{N_{\Lambda}}|_\Lambda=\mathbf{r_\Lambda}/\|\mathbf{r_\Lambda}\|$.

Next, define $w(\xi):=\mathbf{r_\Lambda}(\xi)\cdot\mathbf{N_{\Lambda}}(\xi)$ for $\xi\in\Sigma_\Lambda$, and consider the Cauchy problem:
    \begin{equation}
        \begin{cases}
            \Delta_3 h_\Lambda=0,\\
            h_\Lambda|\Sigma_\Lambda=0,\mathbf{N_{\Lambda}}\cdot\nabla_3 h_\Lambda=w.
        \end{cases}
    \end{equation}
   Here, $\Delta_3:=\partial_{xx}+\partial_{yy}+\partial_{tt}$ denotes the full Laplacian in $\mathbb{R}^3$. The function $w$ is analytic. Since $\Delta_3$ is an elliptic operator, the surface $\Sigma_\Lambda$ is automatically noncharacteristic for the Cauchy problem. Therefore, by the Cauchy-Kovalevskaya Theorem, there exists a unique local solution $h_\Lambda:\Omega\rightarrow\mathbb{R}$ in an open neighborhood $\Omega\subset \mathbb{R}^3$ of~$\Sigma_{\Lambda}$, and in particular, of the initial curve $\Lambda$.

Observe that the function $h_\Lambda$ is zero on $\Sigma_{\Lambda}$, and thus its two tangential derivatives along this surface vanish. Accordingly, at each $\xi\in\Sigma_{\Lambda}$, the gradient $\nabla_3 h_\Lambda$ has the direction of the normal vector $\mathbf{N_{\Lambda}}$. We then conclude from the previous discussion that 
$\nabla_3 h_\Lambda|_\Lambda=\mathbf{r_\Lambda}$, and the lemma follows.

\section{Proof of Lemma~\ref{lemma:FGbeta}}\label{app:appendixB}

Using Equations \eqref{eq:vxomega} and \eqref{eq:vyomega} for $\partial_\sigma v_x$ and $\partial_\sigma v_y$, after tedious but straightforward computations, the reduced normal determinant $\eta$ can be written in terms of the data $\beta$, $F$ and $G$, and on the vectors of the moving frame $\mathcal{B}$, as:
\begin{align}\label{eq:etapol}
    \begin{split}
        \eta(\xi)=&A(\xi)F(\xi)^2+B(\xi)G(\xi)^2+C(\xi)F(\xi)G(\xi)\\[0.8ex]
        &+D(\xi)F(\xi)+E(\xi)G(\xi)+F_0(\xi),\quad\xi\in L_1,
    \end{split}
\end{align}
where the coefficients are:
\begin{align}
    A(\xi) &:= -t_3(\xi)\big(k_1(\xi)^2 + k_2(\xi)^2\big),\\[0.8ex]
    B(\xi) &:= -t_3(\xi)\big(n_1(\xi)^2 + n_2(\xi)^2\big),\\[0.8ex]
    C(\xi) &:= -2t_3(\xi)\big(n_1(\xi)k_1(\xi)+n_2(\xi)k_2(\xi)\big),\\[0.8ex]
    D(\xi) &:= \beta(\xi) t_3(\xi)-2t_3(\xi)\big(X_\beta(\xi)k_1(\xi) + Y_\beta(\xi)k_2(\xi)\big),\label{eq:coeffD}\\[0.8ex]
    E(\xi) &:= -2t_3(\xi)\big(X_\beta(\xi)n_1(\xi) + Y_\beta(\xi)n_2(\xi)\big),\label{eq:coeffE}\\[0.8ex]
    F_0(\xi) &:= \beta(\xi)\big(X_\beta(\xi)n_2(\xi) - Y_\beta(\xi)n_1(\xi)\big)-t_3(\xi)\big(X_\beta(\xi)^2 + Y_\beta(\xi)^2),\quad\xi\in L_1,\label{eq:remaindereta}
\end{align}
 with
\begin{align}
X_\beta(\xi):=t_1k_3 \partial_\theta\beta +\beta\big(t_1\partial_\theta k_3+t_3\partial_\theta k_1+n_3(\partial_\theta\mathbf{k}\times\mathbf{k})_1+2k_3k_1(\mathbf{t}\cdot\partial_\theta\mathbf{k})\big),\quad\xi\in L_1.
\end{align}
and
\begin{align}
Y_\beta(\xi):=t_2k_3 \partial_\theta\beta +\beta\big(t_2\partial_\theta k_3+t_3\partial_\theta k_2+n_3(\partial_\theta\mathbf{k}\times\mathbf{k})_2+2k_3k_2(\mathbf{t}\cdot\partial_\theta\mathbf{k})\big),\quad\xi\in L_1.
\end{align}
All these coefficients, as well as the unknowns $F$ and $G$, are functions of $\xi\in L_1$, so we can interpret Equation \eqref{eq:etapol} as an expression for a parametrized family of second-degree polynomials in $(F,G)$ indexed by $\xi$. To make this rigorous, for all $\xi\in L_1$, we define:
\begin{align}\label{eq:qpoly}
    \begin{split}
    Q_\xi:&\mathbb{R}^2\rightarrow\mathbb{R},\\
    &(a,b)\mapsto A(\xi)a^2+B(\xi)b^2+C(\xi)ab+D(\xi)a+E(\xi)b+F_0(\xi).
    \end{split}
\end{align}
By construction, $\eta(\xi)=Q_{\xi}(F(\xi),G(\xi))$. To prove the lemma, i.e., in order to show that it is possible to choose $F$ and $G$, along with $\beta$, so that $\eta$ is positive on $L_1$, we have to study the global behavior of $Q_{\xi}$.

To this end, we first observe that the determinant of the Hessian matrix of $Q_\xi$ is given by:
\begin{equation}
    \cD_{Q}(\xi):=4t_3(\xi)^4,\quad\xi\in L_1.
\end{equation}
This determinant does not take negative values, and moreover, it is zero only at the horizontal points of $L_1$. If $\xi$ is not horizontal, then $Q_{\xi}$ has a unique (global) extremum $p_c=(a_c,b_c)\in\mathbb{R}^2$. After a few straightforward computations we get the dependence of $p_c$ on $\xi$:
\begin{align}\label{eq:criticalpointQ}
    a_{c}(\xi)&=\frac{1}{2t_3^2}\big(\beta(n_1^2+n_2^2)+2t_3(Y_\beta n_1-X_\beta n_2)\big),\\[0.8ex]
    b_{c}(\xi)&=-\frac{1}{2t_3^2}\big(\beta(k_1n_1+k_2n_2)+2t_3(Y_\beta k_1-X_\beta k_2)\big).
\end{align}
The corresponding critical value as $Q_c(\xi)$ takes the simple form:
\begin{equation}\label{eq:criticalvalueQ}
    Q_c(\xi)=\frac{1}{4 t_3}\beta^2(n_1^2+n_2^2).
\end{equation}

For the construction of $F,G,\beta$ we need to consider the different parts of $L_1$: horizontal points, positive intervals, and negative intervals. Accordingly, let $\xi_j$ be a horizontal point of $L_1$, where $j\in\mathbb{Z}/m\mathbb{Z}$, such that the interval of $L_1$ that goes from $\xi_{j-1}$ to $\xi_j$ is positive and the one that goes from $\xi_j$ to $\xi_{j+1}$ is negative (other cases are analogous). 

First, a direct computation yields:
\begin{equation}
    \eta(\xi_j)=-\beta(\xi_j)\partial_\theta\beta(\xi_j)k_3(\xi_j)^2.
\end{equation}
Since $\Sigma$ is noncharacteristic for the heat equation, the coefficient $k_3(\xi_j)$ is different from zero. Consequently, $\eta(\xi_j)$ is positive if and only if $\beta(\xi_j)$ and $\partial_\theta \beta(\xi_j)$ have opposite nonzero signs. For simplicity, we also choose $F(\xi_j)=G(\xi_j)=0$ for every $j$. By continuity, for each $j$, there exists a small $\tau_j>0$ such that $Q_{\xi}(0,0)>0$ for all $\xi=\phi(\theta)$ with $\theta\in I_j:=[\theta_j-\tau_j,\theta_j+\tau_j]$. Remember that $\theta_j$ denotes the inverse image of $\xi_j$ in $\mathbb{S}^1$.

Second we focus on the interval of $L_1$ that goes from $\xi_{j-1}$ to $\xi_j$. We make $F$ and $G$ take the value zero in a neighborhood of the horizontal points $\xi_{j-1}$ and $\xi_j$. So, for instance we define the subinterval:
\begin{equation}
    U_{j,+}:=\left[\theta_{j-1},\theta_{j-1}+\frac{\tau_{j-1}}{4}\right]\cup\left[\theta_j-\frac{\tau_j}{4},\theta_j\right].
\end{equation}
and we choose:
\begin{equation}\label{eq:zerovaluesFG}
    \begin{cases}
        F(\xi)=0,\\
        G(\xi)=0,
    \end{cases}\quad\forall\xi\in\phi(U_{j,+}).
\end{equation}
Now we work on the complementary interval to $U_{j,+}$:
\begin{equation}
    V_{j,+}:=\left[\theta_{j-1}+\frac{\tau_{j-1}}{4},\theta_j-\frac{\tau_{j}}{4}\right].
\end{equation}
For $\xi\in \phi(V_{j,+})$ the determinant $\cD_Q(\xi)$ is positive, and $A(\xi)$ is negative, so the polynomial $Q_{\xi}$ has a global maximum $p_c$ given by \eqref{eq:criticalpointQ}, whose corresponding critical value $Q_c$ appears in \eqref{eq:criticalvalueQ}. This value is positive   if and only $\beta(\xi)\neq 0$, so we take a function $\beta$ that is nonvanishing for all $\xi$ in the positive interval of $L_1$ that goes from $\xi_{j-1}$ to $\xi_j$. 

To define $F$ and $G$ in $V_{j,+}$, we subdivide this interval into three. The first one is:\begin{equation}
    V_{j,+}^\prime:=\left[\theta_{j-1}+\frac{3\tau_{j-1}}{4},\theta_j-\frac{3\tau_{j}}{4}\right].
\end{equation}
In this interval we set:
\begin{equation}
    \begin{cases}
        F(\xi)=a_c(\xi),\\G(\xi)=b_c(\xi),
    \end{cases}\quad\forall\xi\in \phi(V_{j,+}^\prime).
\end{equation}
The second subinterval is:
\begin{equation}
    V_{j,+}^{\prime\prime}:=\left[\theta_j-\frac{3\tau_{j}}{4},\theta_j-\frac{\tau_{j}}{4}\right].
\end{equation}
Denoting the endpoints of $V_{j,+}^{\prime\prime}$ by
\begin{equation}
    u_{j,+}:=\theta_j-\frac{3\tau_{j}}{4}\quad \text{and} \quad v_{j,+}:=\theta_j-\frac{\tau_{j}}{4}
\end{equation}
we define $F$ and $G$ as:
    \begin{equation}\label{eq:FGgluing}
        \begin{cases}
            F(\xi)=\lambda(\xi)a_c(\xi),\\
            G(\xi)=\lambda(\xi)b_c(\xi),
        \end{cases}
        \quad\text{where}\quad\lambda(\xi):=\left(\frac{\theta(\xi)-u_{j,+}}{v_{j,+}-u_{j,+}}\right),\quad \forall\xi\in\phi(V_{j,+}^{\prime\prime}).
    \end{equation}
We recall that $\theta(\xi)$ denotes the angular coordinate $\theta$ corresponding to the point $\xi$. For each $\xi\in\phi(V_{j,+}^{\prime\prime})$, the point in $\mathbb{R}^2$ with coordinates $(F(\xi),G(\xi))$ is in the straight-line segment between $(0,0)$ and $(a_c(\xi),b_c(\xi))$. Since the polynomial $Q_{\xi}$ is positive at both, and because of its concave behavior for $\xi\in\phi(V_{j,+})$, it follows that $\eta(\xi)$ is indeed positive for this choice of $F$ and $G$ in $V_{j,+}^{\prime\prime}$.

The third interval is:
\begin{equation}
    V_{j,+}^{\prime\prime\prime}:=\left[\theta_{j-1}+\frac{\tau_{j-1}}{4},\theta_{j-1}+\frac{3\tau_{j-1}}{4}\right].
\end{equation}
In this subinterval we do the same type of gluing as previously for $V_{j,+}^{\prime\prime}$, instead that this time we consider a neighborhood of $\xi_{j-1}$ in $L_1$. The expression for $F$ and $G$ here is analogous to the one appearing in \eqref{eq:FGgluing}, so we give no further details.

In summary, taking into account the previous choice of $\beta$, we have constructed continuous functions $F$ and $G$, on the positive interval of $L_1$, such that $\eta(\xi)$ is positive in the closure of such an interval.

To complete the construction it remains to consider the interval of $L_1$ that goes from $\xi_j$ to $\xi_{j+1}$ (the negative interval). The construction here mirrors that for the positive interval, so we only sketch it. First, we set:
\begin{equation}
    U_{j,-}:=\left[\theta_{j},\theta_{j}+\frac{\tau_{j}}{4}\right]\cup\left[\theta_{j+1}-\frac{\tau_{j+1}}{4},\theta_{j+1}\right],
\end{equation}
and we choose:
\begin{equation}
    \begin{cases}
        F(\xi)=0,\\
        G(\xi)=0,
    \end{cases}
    \quad\forall\xi\in\phi(U_{j,-}).
\end{equation}

Next, introducing the complementary interval to $U_{j,-}$: 
\begin{equation}
    V_{j,-}:=\left[\theta_j+\frac{\tau_j}{4},\theta_{j+1}-\frac{\tau_{j+1}}{4}\right],
\end{equation}
 we easily infer that for $\xi\in \phi(V_{j,-})$, the polynomial $Q_{\xi}$ is strictly convex. The critical point $p_c(\xi)$ is a global minimum, and the corresponding critical value $Q_c(\xi)$ is always nonpositive. Contrary to the case of the the positive interval, the function $\beta$ is allowed to vanish here.

For each $\xi\in\phi(V_{j,-})$, define the set
\begin{equation}
    E_j(\xi):=\big\{(a,b)\in\mathbb{R}^2:Q_{\xi}(a,b)\leq 0\big\},   
\end{equation}
which is a nonempty compact subset of $\mathbb R^2$, and $E:=\bigcup_{\xi\in\phi( V_{j,-})}E_j(\xi)$, which is also a compact subset of $\mathbb R^2$ because $V_{j,-}$ is compact. As before, take:
\begin{equation}
    V_{j,-}^\prime:=\left[\theta_j+\frac{3\tau_j}{4},\theta_{j+1}-\frac{3\tau_{j+1}}{4}\right],
\end{equation}
and a point $(M_1,M_2)\notin E$. Then, we choose:
\begin{equation}
    \begin{cases}
        F(\xi)=M_1,\\
        G(\xi)=M_2,
    \end{cases}
    \quad\forall\xi\in\phi(V_{j,-}^\prime).
\end{equation}
Similarly,
\begin{equation}
    V_{j,-}^{\prime\prime}:=\left[\theta_j+\frac{\tau_j}{4},\theta_j+\frac{3\tau_j}{4}\right].
\end{equation}
In this subinterval, we glue the values of $F$ and $G$ in $U_{j,-}$ with those we just described in $V_{j,-}^{\prime}$. Here we cannot just use a linear interpolation since the region where $Q_{\xi}$ is positive is no longer convex, but we can easily find a continuous path $\mathbb R^2\backslash E$ that connects $(0,0)$ and $(M_1,M_2)$. With it, we can easily define $F(\xi)$ and $G(\xi)$ for $\xi\in V_{j,-}^{\prime\prime}$ so that $Q_\xi(F(\xi),G(\xi))>0$ in such a subinterval.

Finally, in the subinterval
\begin{equation}
     V_{j,-}^{\prime\prime\prime}:=\left[\theta_{j+1}-\frac{3\tau_{j+1}}{4},\theta_{j+1}-\frac{\tau_{j+1}}{4}\right].
\end{equation}
 we define functions $F,G$ as in the case of the subinterval $V_{j,-}^{\prime\prime}$.

Summarizing, the construction above yields two $C^0$ functions $F$ and $G$ that give a positive reduced normal determinant $\eta$, provided that $\beta(\xi)$ is different from zero for all $\xi\in L_1^+$, and that itself and its derivative $\partial_\theta\beta$ have opposite nonzero signs at the horizontal points. Furthermore, taking into account Constraint \textbf{(B)} and the decomposition of $L_1$ in Section \ref{section:surface}, it is required that $\beta(\xi)$ is negative for $\xi\in L_{1,\max}^+$ and positive for $\xi\in L_{1,\min}^+$. An analytic function $\beta:\Omega\rightarrow\mathbb{R}$ such that $\beta|_{L_1}$ satisfies all these conditions does clearly exist. Concerning $F$ and $G$, using {Whitney's Analytic Approximation Theorem} \cite[Theorem 5.6.]{hirsch_differential_1976}, we can take two analytic functions that are arbitrarily close to the previous continuous functions $F$ and $G$, for which $\eta(\xi)$ is still positive for all $\xi\in L_1$. This completes the proof of the lemma.

\bibliographystyle{amsplain}

\end{document}